\documentclass{amsart}
\usepackage[english]{babel}
\usepackage[latin1]{inputenc}
\usepackage{graphicx}
\usepackage{amsmath}
\usepackage{amsfonts}
\usepackage{amssymb}
\usepackage{amsthm}
\usepackage{hyperref}

\newtheorem{thm}{Theorem}[section]

\newtheorem{lem}[thm]{Lemma}
\newtheorem{prop}[thm]{Proposition}
\theoremstyle{definition}
\newtheorem{defn}[thm]{Definition}
\theoremstyle{remark}
\newtheorem{rem}[thm]{Remark}

\numberwithin{equation}{section}

\begin{document}

\title{Partial Degree Formulae for Plane Offset Curves}%
\author{F. San Segundo}
\author{J. R. Sendra}
\thanks{Both authors supported by the Spanish `` Ministerio de Educaci\'on y
Ciencia" under the Project MTM2005-08690-C02-01 and  by the ``Direcci\'on General de
Universidades de la Consejer\'{\i}a de Educaci\'on de la CAM y la Universidad de Alcal\'a"
under the project CAM-UAH2005/053.}
\address{Departamento de Matem\'aticas, Universidad de Alcal\'a, E-28871-Madrid, Spain}
\begin{abstract}
In this paper we present several formulae for computing the partial degrees of the
defining polynomial of the offset curve to an irreducible affine plane curve
 given implicitly, and we see how these formulae
particularize to the case of rational curves. In addition, we present a formula for
computing  the degree w.r.t the distance variable.
\end{abstract}
\email{fernando.sansegundo@uah.es}%
\maketitle

% main text
%%%%%%%%%%%%%%%%%%%%%%%%%%%%%%%%%%%%%%%%%%%%%%%%%%%%%%%%%%%%%%%%%%%%%%%%%%%%%%%%%%%%%%%%%
%%%%%%%%%%%%%%%%%%%%%%%%%%%%%%%%%%%%%%%%%%%%%%%%%%%%%%%%%%%%%%%%%%%%%%%%%%%%%%%%%%%%%%%%%
\section{Introduction}
\label{SecIntroduction}
%%%%%%%%%%%%%%%%%%%%%%%%%%%%%%%%%%%%%%%%%%%%%%%%%%%%%%%%%%%%%%%%%%%%%%%%%%%%%%%%%%%%%%%%%
%%%%%%%%%%%%%%%%%%%%%%%%%%%%%%%%%%%%%%%%%%%%%%%%%%%%%%%%%%%%%%%%%%%%%%%%%%%%%%%%%%%%%%%%%
Offset curves and surfaces are well-known geometric objects in the field of computer
aided geometric design, possibly  because they constitute a powerful tool in many
applications (see \cite{handbook}, \cite{HL97}, \cite{Pet98}). On the other hand, offset
construction is a real mathematical challenge. Even though one starts from a very simple
curve/surface, the offset is usually much more complicated. Because of this fact, many
authors try to deduce a priori information (on applied, algorithmic, or even theoretical
aspects) of the offset from the original generating curve/surface. For instance,
relevant results have been achieved in problems like:  the determination of the genus of
the offset (see \cite{ASS97}), deciding the rationality and parametrizing offsets (see
\cite{ASS96}, \cite{L95a}, \cite{Pet98}, \cite{Po95}, \cite{SS00}), implicitization
techniques (see \cite{Hof90}, \cite{Hof93}, \cite{Wang}), analyzing its topological type
(see \cite{AS06}, \cite{FN90a} ), studying analytic and algebraic properties (see
\cite{FN90a}, \cite{FN90b}, \cite{SS99}), etc.

An additional problem, not mentioned above and that is the central topic of this paper,
is the computation of the degree of the offset. Results in this direction, for offset
curves, can be found in \cite{FN90b} for the parametric case, in \cite{SSS05} for the
implicit and parametric case, and in \cite{Anton} for the implicit case. Note that the
knowledge of the offset degree can be applied, for instance, for constructing ad hoc
offset implicitization algorithms based on interpolation techniques.

All the contributions mentioned above deal with the problem of computing the degree of
the offset curve; that is the total degree of its defining polynomial. In this paper, we
complete this analysis providing formulae for the partial degree of the offset defining
polynomial w.r.t. each variable, including the distance one. This extension of the work
presented in \cite{SSS05} may have relevant implications in the improvement of
interpolation-based algorithms for implicitizing, since with these additional
information the interpolation space is reduced.

In order to formally state the problem, we consider a polynomial $g(x_1,x_2,d)$ in the
variables $\{x_1,x_2,d\}$ such that for all values $d_0$ of $d$, but either none or
finitely many exceptions, $g(x_1,x_2,d_0)$ is the implicit equation of the offset at
distance $d_0$; this polynomial is called the {\sf generic offset equation} and its
existence and specialization properties are established in Section
\ref{SecOffsetEquation}. In this situation, the problem consists in computing the
partial degrees $\deg_{x_1}(g)$, $\deg_{x_2}(g)$, and $\deg_d(g)$. Concerning to the
coordinate partial degrees, i.e. $\deg_{x_1}(g)$, $\deg_{x_2}(g)$, we present four
different formulae; two of them for the implicit cases, and the two others for the
parametric case. The distance degree formula is stated assuming that the input generator
curve is given by means of its implicit equation.

The strategy we follow for developing the formulae is essentially the one used in
\cite{SSS05}. That is, we consider the intersection of the offset with a general
vertical/horizonal line. Then, the partial degree is the number of intersection points.
This number of intersection points is deduced from  the intersection points of the
original curve with an auxiliary curve, directly deduced from the input, and constructed
ad hoc for each degree problem. Therefore, explicit knowledge on the offset is avoided.
Note that the main difference, of the reasoning here and the reasoning in \cite{SSS05},
is that the total degree of a curve is the number of intersections with a generic line
but, for the partial degrees, generic vertical or horizontal lines need to be
considered.

The structure of the paper is as follows. In Section \ref{SecOffsetEquation} we
introduce the notion of generic offset and generic offset equation, and we establish
their main properties. In Section \ref{SecStrategyPartial} we describe the theoretical
strategy for computing the partial degree formulae. In Section \ref{secCurvaS} we
introduce the auxiliary curve $\mathcal S$ as well as the fake and non-fake intersection
points. Finally, in Section \ref{secPartialDegreeFormulae} we apply these ideas to
develop the partial degree formulae for the implicit case. The particularization of
these formulae to the parametric case is done in Section \ref{secParametricCase} After
that, the paper focuses on the distance degree formula. This is done in two sections. In
Section \ref{SecStrategyDistance} we show how to adapt the strategy for this special
case, and in Section \ref{secDistanceFormula} the distance degree formula is deduced.
The papers ends with an appendix (in page \pageref{Apenddix}) where all the degrees
(total and partial) are listed for a collection of curves.

%%%%%%%%%%%%%%%%%%%%%%%%%%%%%%%%%%%%%%%%%%%%%%%%%%%%%%%%%%%%%%%%%%%%%%%%%%%%%%%%%%%%%%%%%
%%%%%%%%%%%%%%%%%%%%%%%%%%%%%%%%%%%%%%%%%%%%%%%%%%%%%%%%%%%%%%%%%%%%%%%%%%%%%%%%%%%%%%%%%
\section{The generic equation of the offset}
\label{SecOffsetEquation}
%%%%%%%%%%%%%%%%%%%%%%%%%%%%%%%%%%%%%%%%%%%%%%%%%%%%%%%%%%%%%%%%%%%%%%%%%%%%%%%%%%%%%%%%%
%%%%%%%%%%%%%%%%%%%%%%%%%%%%%%%%%%%%%%%%%%%%%%%%%%%%%%%%%%%%%%%%%%%%%%%%%%%%%%%%%%%%%%%%%

We start recalling the classical and intuitive concept of offset curve. This notion will
be formalized in this section. Let ${\mathcal C}$ be a plane curve, and let
$p\in{\mathcal C}$. Let ${\mathcal L}_N$ be the normal line to ${\mathcal C}$ at $p$
(assume for now that this normal line is well defined). Let $q_1, q_2$ be the two points
of ${\mathcal L}_N$ at a fixed distance $d_0\in{\mathbb C}^*$ of $p$.  Then, the offset
curve (or parallel curve) to ${\mathcal C}$ at distance $d_0$, is the set ${\mathcal
O}_{d_0}({\mathcal C})$ of the points $q_i$ obtained by means of this geometric
construction.

As the distance $d_0$ varies, different offset curves are obtained. The idea is to have
a global expression of the offset for all (or almost all) distances. This motivates the
concept of {\sf generic equation of the offset to $\mathcal C$}. This generic equation
is a polynomial, depending on the variable distance $d$, such that for every (or almost
every, see the examples below) value of $d$, the equation specializes to the equation of
the offset at that particular distance.

Using this informal definition of generic offset equation, and using Gr\"{o}bner basis
techniques, one can see that if $\mathcal C$ is the parabola $y_2-y_1^2=0$, then the
generic equation of its offset is: \vspace{1mm}

 \noindent
 $g(x_1,x_2,d)=-48\,{d}^{2}{x_1}^{4}-32\,{d}^{2}{x_1}^{2}{x_2}^{2}+48
\,{d}^{4}{x_1}^{2}+16\,{x_1}^{6}+16\,{x_2}^{2}{x_1 }^{4}+16\,{d}^{4}{x_2}^{2}
-16\,{d}^{6}-40\,x_2\,{x_1}^{ 4}-32\,{x_1}^{2}{x_2}^{3}+8\,{d}^{2}x_2\,{x_1}^{2}
-32\,{d}^{2}{x_2}^{3}+32\,{d}^{4}x_2+{x_1}^{4}+32\,{{ \it
x1}}^{2}{x_2}^{2}+16\,{x_2}^{4}-20\,{d}^{2}{x_1}^{2
}-8\,{d}^{2}{x_2}^{2}-8\,{d}^{4}-2\,x_2\,{x_1}^{2}-8\,{
x_2}^{3}+8\,x_2\,{d}^{2}+{x_2}^{2}-{d}^{2}.$
 \vspace{1mm}

\noindent In addition, and using again Gr\"{o}bner basis techniques, one may check that for
every distance the generic offset equation specializes properly. However, the generic
offset equation of  the circle $y_1^2+y_2^2-1=0$ factors as the product of two circles
of radius $1+d$ and $1-d$, that is:
\[g(x_1,x_2,d)= \left( x_1^2+x_2^2-(1+d)^2\right) \left( x_1^2+x_2^2-(1-d)^2\right).
\]
Observe that for $d_0=1$, this generic equation gives
\[g(x_1,x_2,1)= \left( x_1^2+x_2^2-2^2\right) \left( x_1^2+x_2^2\right)=
\left( x_1^2+x_2^2-2^2\right)\left( x_1+i x_2\right)\left( x_1-i x_2\right)
\]
which describes the union of a circle of radius $2$, and two complex lines. This is not
a correct representation of the offset at distance $1$ to $\mathcal C$, which consists
of the union of the circle of radius $2$ and a point (the origin). Thus, in this example
we see that the generic offset equation does not specialize properly for $d_0=1$.
Nevertheless, for every other value of $d_0$ the specialization is correct.

In these examples we have introduced some of the {\bf notation} that we will use in the
sequel. The variables $\bar y=(y_1,y_2)$ will be used for the equation of the curve
$\mathcal C$, and $\bar x=(x_1,x_2)$ will be used for the equation of the offset to
$\mathcal C$, both for a particular distance or generically. The implicit equation of
$\mathcal C$ is $f(y_1,y_2)=0$ and the generic offset equation is $g(x_1,x_2,d)=0$.

After these examples, we proceed to formally introduce the notions of offset and of
generic offset equation. This can be done using a geometrical approach, by means of
incidence diagrams (see \cite{SS99}), or equivalently using results from Elimination
Theory. Here we follow this second approach. For this purpose, let $\mathcal C$ be an
irreducible algebraic plane curve given by the polynomial $f(y_1,y_2)\in{\mathbb
C}[y_1,y_2]$ such that $f$ does not divide to $f_{1}^{2}+f_{2}^{2}$. Note that this
implies that the set of non-isotropic points of $\mathcal C$ is open and non-empty (see
Proposition 2 in \cite{SS99}); i.e. the set of points of $\mathcal C$ at which the
non-zero normal vectors $(n_1,n_2)$ satisfies that $n_{1}^{2}+n_{2}^{2}\neq 0$.
Moreover, by Proposition 1 in \cite{SS99}, if $\mathcal C$ is real and irreducible this
condition holds. Consider the following polynomial system:
\[
\left.\begin{array}{lr}
&f(y_1,y_2)=0\\
b(\bar y, \bar x,d):& (x_1-y_1)^2+(x_2-y_2)^2-d^2=0\\
n(\bar y, \bar x):&-f_2(\bar y)(x_1-y_1)+f_1(\bar y)(x_2-y_2)=0\\
w(\bar y,u):& u\cdot(f_1^2(\bar y)+f_2^2(\bar y))-1=0
\end{array}\right\}
\equiv\mathfrak{S}_1(d)
\]
where $f,b,n,w\in{\mathbb C}[\bar y, \bar x,d,u]$, with $\bar x=(x_1,x_2), \bar
y=(y_1,y_2)$ and $f_i=\dfrac{\partial f}{\partial y_i}$.

Note that $d$ is considered here as a variable, representing the distance. The second
equation, $b(\bar y, \bar x,d)$, represents a circle of radius $d$ centered at the point
$\bar y\in{\mathcal C}$, and the third one defines the normal line to $\mathcal C$ at
$\bar y$. The last equation excludes the possibility of $\bar y$ being a singular (or,
in general, isotropic) point of $\mathcal C$. In addition, observe that we have assumed
that $f$ does not divide to $f_{1}^{2}+f_{2}^{2}$, and therefore ${\mathfrak S_1}(d)$
has always solutions.

First, we will establish the existence of the generic equation of the offset.  Let
\[I(d)=<f(\bar y),b(\bar y, \bar x,d),n(\bar y, \bar x),w(\bar y,u)>\]
be the ideal in ${\mathbb C}[\bar y, \bar x,d,u]$ generated by the polynomials $\{f, b,
n, w\}$. We denote by
\[\Omega(d)=\mathbf V(I(d))\subset{\mathbb C}^6\]
the affine algebraic set defined by $I(d)$; that is, $\Omega(d)$ is the set of solutions
in $\mathbb C^6$ of the system $\mathfrak{S}_1(d)$.

Now, for every particular $d_0\in{\mathbb C}^*$, let
\[I(d_0)=<f(\bar y),b(\bar y, \bar x,d_0),n(\bar y, \bar x),w(\bar y,u)>\]
be the ideal in ${\mathbb C}[\bar y, \bar x,u]$ generated by $\{f, b(d_0), n, w\}$.  And
let
\[\Omega(d_0)=\mathbf V(I(d_0))\subset{\mathbb C}^5\]
be the affine algebraic set defined by $I(d_0)$.

We consider the following two projection maps:
\begin{enumerate}
\item[] $\pi:{\mathbb C}^6\rightarrow{\mathbb C}^3;(\bar y, \bar x,d,u)\mapsto(\bar x,d)$\hfill {\sf (non-specialized projection)}
\item[] $\pi_0:{\mathbb C}^5\rightarrow{\mathbb C}^2;(\bar y, \bar x,u)\mapsto \bar x$\hfill {\sf (specialized projection)}
\end{enumerate}
In this situation, if one denotes by ${\mathcal A}^*$ the Zariski closure of a set
${\mathcal A}$, one has the following definition:
\begin{defn}
The {\sf offset} to the curve $\mathcal C$ at a distance $d_0$ is
\[{\mathcal O}_{d_0}(\mathcal C)=\left(\pi_0\left(\Omega(d_0)\right)\right)^*\subset{\mathbb C}^2\]
The {\sf generic offset} to the curve $\mathcal C$ is
\[{\mathcal O}_{d}(\mathcal C)=\left(\pi\left(\Omega(d)\right)\right)^*\subset{\mathbb C}^3\]
\end{defn}
\begin{rem}\label{remarkIdealEliminacion}
Note that this means that
\[{\mathcal O}_{d}(\mathcal C)=\mathbf V(\tilde I(d))\]
where $\tilde I(d)=I(d)\cap{\mathbb C}[\bar x, d]$ is the $(\bar y,u)$-elimination ideal
of $I(d)$. Similarly
\[{\mathcal O}_{d_0}(\mathcal C)=\mathbf V(\tilde I(d_0))\]
where $\tilde I(d_0)=I(d_0)\cap{\mathbb C}[\bar x]$ (see \cite{Cox}, Closure Theorem, p.
122).
\end{rem}

The following result guarantees the existence of an equation for the generic offset.
\begin{lem}\label{Lemma01}
${\mathcal O}_{d}(\mathcal C)$ is a surface in ${\mathbb C}^3$.
\end{lem}
{\bf Proof.} {\small This proof follows the reasoning of the proof of Lemma 1 in
\cite{SS99}. Let $K$ be a component of $\Omega(d)$, and let $(p,q,u_0,d_0)\in K$. Since
$w(p,u_0)$=0, $p\in\mathcal C$ is non-isotropic. Moreover, $q\in{\mathcal
O}_{d_0}(\mathcal C)$. Take $P(t)=(x(t),y(t))$ to be a place of $\mathcal C$ centered at
$p$ ($P(t)$ is a local parametrization of $\mathcal C$ by power series). Let $N(t)$ be
the associated normal vector, and let $Q(t)$ be the lifting of $P(t)$ to $q\in{\mathcal
O}_{d_0}(\mathcal C)$ whose center is $q$. That is,
\[Q(t)=P(t)\pm d\dfrac{N(t)}{\|N(t)\|}\]
The choice of sign is decided with the condition that $Q(t)$ is centered at $q$.
Moreover, note that since $p$ is non-isotropic, then $Q(t)$ is also a local
parametrization by power series. Then
\[R(t,d)=\left(P(t),Q(t),d,\dfrac{1}{\|N(t)\|^2}\right)\]
is a local parametrization of $K$ at $(p,q,u_0,d_0)$. It follows that $\dim
K=2$.\hfill\qed}

Therefore ${\mathcal O}_{d}(\mathcal C)$ is defined by a polynomial in ${\mathbb C}[\bar
x,d]$ (see \cite{Shafarevich77}, p.69, Th.3). Thus, we arrive at the following
definition:
\begin{defn}
The {\sf generic offset equation} is the defining polynomial of the surface ${\mathcal
O}_{d}(\mathcal C)$. In the sequel, {\sf we denote} by  $g(x_1,x_2,d)=0$ the generic
offset equation.
\end{defn}

\begin{rem}\label{Remarkg}
\begin{enumerate}
\item[]
\item Observe that the polynomial $g$ may be reducible (recall the example of the circle)
but by construction it is always square-free. Moreover, $g$ is either irreducible or
factors into two irreducible factors not depending only on $d$; this is so because,
generically in $d$, the offset has at most two irreducible components (see \cite{SS99},
Theorem 1).
\item It might happen that $g(\bar x,d)$ has a factor in $\mathbb C[d]$. In order to
avoid this, and w.l.o.g., we will take the generic offset equation to be primitve w.r.t.
$\bar x$.
\end{enumerate}
\end{rem}

The following theorem gives the fundamental property of the generic offset.
\begin{thm}\label{Prop01}
For all but finitely many exceptions, the generic offset equation specializes properly.
That is, there exists a finite (possibly empty) set $\Upsilon\subset{\mathbb C}$ such
that if $d_0\not\in\Upsilon$, then
\[g(x_1,x_2,d_0)=0\]
is the equation of ${\mathcal O}_{d_0}(\mathcal C)$.
\end{thm}

{\bf Proof.}{\small Since $g(\bar x,d)$ defines the equation of ${\mathcal
O}_{d}(\mathcal C)$, and
\[{\mathcal O}_{d}(\mathcal C)=\mathbf V(\tilde I(d))\]
where $\tilde I(d)=I(d)\cap{\mathbb C}[\bar x, d]$ is the $(\bar y,u)$-elimination ideal
of $I(d)$ (see Remark \ref{remarkIdealEliminacion}), it follows that if $G(d)$ is  a
Gr\"{o}bner basis of $I(d)$ w.r.t. an elimination ordering that eliminates $(\bar y,u)$,
then up to multiplication by a non-zero constant, $G(d)\cap{\mathbb C}[\bar x,
d]=\{g(\bar x, d)\}$ is a Gr\"{o}bner basis of $\tilde I(d)$. But then (see \cite{Cox},
exercise 7, page 284) there is a finite (possibly empty) set $\Upsilon\subset\mathbb C$
such that for $d_0\not\in\Upsilon$, $G(d_0)$ specializes well to a Gr\"{o}bner basis of
$I(d_0)$. It follows that, since $\tilde I(d_0)=I(d_0)\cap{\mathbb C}[\bar x]$, then
$G(d_0)\cap{\mathbb C}[\bar x]=\{g(\bar x, d_0)\}$ is a Gr\"{o}bner basis of $\tilde
I(d_0)$. Thus, for $d_0\not\in\Upsilon$, $g(\bar x, d_0)$ is the equation of ${\mathcal
O}_{d_0}(\mathcal C)=\mathbf V(\tilde I(d_0))$.\hfill\qed}

\begin{rem}\label{remExtensionHypersurfaces}
Note that all the results in this section, though they have been presented for plane
curves, extend naturally to the case of offsets to irreducible hypersurfaces (over
algebraically closed fields of characteristic zero).
\end{rem}

%%%%%%%%%%%%%%%%%%%%%%%%%%%%%%%%%%%%%%%%%%%%%%%%%%%%%%%%%%%%%%%%%%%%%%%%%%%%%%%%%%%%%%%%%
%%%%%%%%%%%%%%%%%%%%%%%%%%%%%%%%%%%%%%%%%%%%%%%%%%%%%%%%%%%%%%%%%%%%%%%%%%%%%%%%%%%%%%%%%
\section{Strategy description for the partial degree formulae}\label{SecStrategyPartial}
%%%%%%%%%%%%%%%%%%%%%%%%%%%%%%%%%%%%%%%%%%%%%%%%%%%%%%%%%%%%%%%%%%%%%%%%%%%%%%%%%%%%%%%%%
%%%%%%%%%%%%%%%%%%%%%%%%%%%%%%%%%%%%%%%%%%%%%%%%%%%%%%%%%%%%%%%%%%%%%%%%%%%%%%%%%%%%%%%%%

First we deal with the problem of computing the partial degree in $x_i$ of the generic
offset equation $g(\bar x,d)$. Let $\delta_i$ be the partial degree in $x_i$ of $g$. We
will describe how to compute $\delta_1$. Then, simply exchanging the variables $x_1$ and
$x_2$ allows to compute $\delta_2.$ Also, we will exclude w.l.o.g. in our analysis the
case where $\mathcal C$ is a line. Note that, in particular, this implies that
$\delta_i>0$ in all cases.

When analyzing the offset total degree problem in our previous paper \cite{SSS05}, the
basic idea was to indirectly determine the number of intersection points between a
generic line and the offset ${\mathcal O}_d({\mathcal C})$. Here, for the partial degree
problem, we follow a similar strategy. However, in order to compute $\delta_1$, the
generic line must be horizontal. Let therefore
$$\ell(\bar x,k): x_2-k=0$$
be the equation of a generic horizontal line ${\mathcal L}(k)$. Since the generic offset
equation is not known, we compute indirectly the number of points in ${\mathcal
O}_d({\mathcal C})\cap{\mathcal L}(k)$, by counting the points in ${\mathcal C}$ that,
in a 1:1 correspondence, generate the points in ${\mathcal O}_d({\mathcal
C})\cap{\mathcal L}(k)$. For this purpose, we analyze the solutions of system
$\mathfrak{S}_1(d)$ lying on the line ${\mathcal L}(k)$. That is, the solutions of the
system:
\[
\left.\begin{array}{lr}
&f(y_1,y_2)=0\\
b(\bar y, \bar x,d):& (x_1-y_1)^2+(x_2-y_2)^2-d^2=0\\
n(\bar y, \bar x):&-f_2(\bar y)(x_1-y_1)+f_1(\bar y)(x_2-y_2)=0\\
w(\bar y,u):& u\cdot(f_1^2(\bar y)+f_2^2(\bar y))-1=0\\
\ell(\bar x,k):& x_2-k=0
\end{array}\right\}
\equiv\mathfrak{S}_2(d,k)
\]
The following result provides the theoretical foundation of our strategy, by
establishing the 1:1 correspondence between the points in ${\mathcal O}_d({\mathcal
C})\cap{\mathcal L}(k)$, and the points in ${\mathcal C}$ that generate them.

We recall that a {\sf ramification point} of a curve is a point on the curve where at
least one of the partial derivatives of the implicit equation vanishes. In our case,
since we are analyzing the partial degree $\delta_1$, by abuse of notation, whenever we
speak about ramification points we mean a ramification point where the partial
derivative w.r.t. $y_2$ vanishes.
\begin{thm}\label{Theorem01}
There exists a non-empty Zariski open subset $\Delta$ of ${\mathbb C}^2$ such that for
$(d_0,k_0)\in\Delta$:
\begin{enumerate}
\item There exist exactly $\delta_1$ solutions $\Gamma=\{(p_i,q_i,u_i)\}_{i=1,\ldots,\delta_1}$
of $\mathfrak{S}_2(d_0,k_0)$ satisfying that:
\begin{enumerate}
\item $q_1,\ldots,q_{\delta_1}$ are all different and
${\mathcal L}(k_0)\cap{\mathcal O}_{d_0}({\mathcal C})=\{q_1,\ldots,q_{\delta_1}\}$.
\item $p_1,\ldots,p_{\delta_1}$ are different regular
non-ramification points  of ${\mathcal C}$.
\end{enumerate}
\item None of the points in ${\mathcal C}\cap{\mathcal L}(k_0)$
is a ramification point of ${\mathcal C}$.
\end{enumerate}
\end{thm}
{\bf Proof. }{\small Let us consider the generic offset equation as a polynomial in
$\mathbb C[x_2,d][x_1]$, by writing:
\[g(x_1,x_2,d)=\sum_{i=0}^{\delta_1}g_i(x_2,d)x_1^i,\]
where $g_{\delta_1}$ is not identically zero. Observe that by assumption $\delta_1>0$.
Thus, the set of solutions of $g_{\delta_1}(k,d)=0$ is either empty, or a curve $\Psi_1$
in ${\mathbb C}^2$. We define $\Delta_1={\mathbb C}^2\setminus\Psi_1$.

Besides, by Theorem \ref{Prop01}, we know that there is only a finite set of {\em bad}
distances, $\Upsilon=\{d_1,\ldots,d_m\}$, such that for $d_0\not\in\Upsilon$, the
equation of ${\mathcal O}_{d_0}({\mathcal C})$ is $g(x_1,x_2,d_0)=0$. Let $\Psi_2$ be
the union of the lines with equations $d=d_i$ for $d_i\in\Upsilon$. We define
$\Delta_2=\Delta_1\setminus\Psi_2$. Then, for $(d_0,k_0)\in\Delta_2$,
$$g(x_1,d_0,k_0)=\sum_{i=1}^{\delta_1}g_i(d_0,k_0)x_1^i=0$$
is a polynomial in $x_1$ of degree $\delta_1$ (the leading coefficient does not vanish
because of the construction of $\Delta_1$). Now, since $g$ is square-free (see Remark
\ref{Remarkg}), $\operatorname{Dis}_{x_1}(g(x_1,k,d))$ is a non-identically zero
polynomial in $(k,d)$. Thus, it defines a curve $\Psi_3$ in the $(k,d)-$plane. We define
$\Delta_3=\Delta_2\setminus\Psi_3$.

Let now $\sigma=(\sigma_1,\sigma_2)$ be one of the finitely many singularities or
vertical ramification points of ${\mathcal C}$ (that is, one of the finitely many
solutions of $f=f_2=0$; note that ${\mathcal C}$ is irreducible). We compute the
following resultant between the generic offset polynomial and the equation of a
$d$-circle centered at $\sigma$.
\[R_\sigma(k,d)=\operatorname{Res}_{x_1}(g(x_1,k,d),(x_1-\sigma_1)^2+(k-\sigma_2)^2-d^2)\]
This resultant can only vanish identically if both polynomials have a common factor in
$x_1$. But the polynomial defining the circle is irreducible. Thus, this could only
happen if, for every $d_0\not\in\Upsilon$, ${\mathcal O}_{d_0}({\mathcal C})$ contains a
circle of radius $d_0$ centered at $\sigma$. This would imply that ${\mathcal C}$ is
itself a circle centered at $\sigma$, which is impossible since $\sigma\in{\mathcal C}$.
Thus, $R_{\sigma}$ is not zero, and it defines a curve in  ${\mathbb C}^2$. Let $\Psi_4$
be the curve obtained as the union of such curves for all the possible points $\sigma$.
We define $\Delta_4=\Delta_3\setminus\Psi_4$.

Now, observe that for $(d_0,k_0)\in\Delta_4$, no intersection point of ${\mathcal
O}_{d_0}({\mathcal C})$ and ${\mathcal L}(k_0)$ can be associated with a singularity or
vertical ramification point of ${\mathcal C}$.

Since ${\mathcal C}$ has only finitely many vertical ramification points, we can exclude
those values of $k$ such that the line $x_2=k$ passes through one of those vertical
ramification points. Let $\Psi_5$ be the finite union of such lines, and define
$\Delta_5=\Delta_4\setminus\Psi_5$.

Take
\[\Delta=\Delta_5\]
Then, if $(d_0,k_0)\in\Delta$, because of the construction of $\Delta_2$, we know that
$g(x_1,x_2,d_0)$ is the equation of ${\mathcal O}_{d_0}({\mathcal C})$. Besides, the
equation
$$g(x_1,d_0,k_0)=\sum_{i=1}^{\delta_1}g_i(d_0,k_0)x_1^i=0$$
has exactly $\delta_1$ different roots because of the construction of $\Delta_1$ and
$\Delta_3$. Every solution of this equation represents an affine intersection point of
${\mathcal O}_{d_0}({\mathcal C})$ and $\mathcal L(k_0)$. Moreover, because of the
choice of $\Delta_4$, these points are associated to regular non-ramification affine
points of ${\mathcal C}$. This proves statement (1) of the theorem. Moreover, for
$(d_0,k_0)\in\Delta$ the system $f(\bar y)=0, f_2(\bar y)=0, y_2=k_0$ has no solutions,
because of the construction of $\Delta_5$. This proves statement (2). \hfill\qed}

\begin{rem}
\begin{enumerate}
\item[]
\item In the sequel we {\bf assume} that for $(d_0,k_0)\in\Delta$, $g(\bar y,d_0)=0$ is the implicit
equation of ${\mathcal O}_{d_0}({\mathcal C})$. This can be assumed w.l.o.g., simply
replacing $\Delta$ by $\Delta\setminus \left[({\mathbb
C}\setminus\Upsilon)\times{\mathbb C}\right]$ (see Theorem \ref{Prop01}).
\item  Note that besides the $\delta_1$ solutions mentioned in the theorem, the system
$\mathfrak{S}_2(d_0,k_0)$ may have other solutions. We will analyze in the next section
the distinction between these two types of solutions of the system.
\end{enumerate}
\end{rem}

We have seen that, generically in $k$ and $d$, every point $q_j\in{\mathcal
O}_d({\mathcal C})\cap{\mathcal L}(k)$ is associated to a regular affine point $p_j\in
{\mathcal C}$, and this correspondence is a bijection. The number of such points is the
offset partial degree $\delta_1$. The strategy now is to eliminate $x_1,x_2$ from the
system $\mathfrak{S}_2(d,k)$ in order to obtain information about $\delta_1$ through the
solutions $(y_1,y_2)$ of the resulting system. This means that we switch our attention
from the points $q=(x_1,x_2)\in{\mathcal O}_d\cap{\mathcal L}(k)$ to the associated
points $p=(y_1,y_2)\in{\mathcal C}$. In order to do that we will identify these
associated points as intersection points of ${\mathcal C}$ with a certain auxiliary
curve ${\mathcal S}$ (see Definition \ref{defCurvaS} below).

%%%%%%%%%%%%%%%%%%%%%%%%%%%%%%%%%%%%%%%%%%%%%%%%%%%%%%%%%%%%%%%%%%%%%%%%%%%%%%%%%%%%%%%%%
%%%%%%%%%%%%%%%%%%%%%%%%%%%%%%%%%%%%%%%%%%%%%%%%%%%%%%%%%%%%%%%%%%%%%%%%%%%%%%%%%%%%%%%%%
\section{The Auxiliary Curve ${\mathcal S}$}\label{secCurvaS}
%%%%%%%%%%%%%%%%%%%%%%%%%%%%%%%%%%%%%%%%%%%%%%%%%%%%%%%%%%%%%%%%%%%%%%%%%%%%%%%%%%%%%%%%%
%%%%%%%%%%%%%%%%%%%%%%%%%%%%%%%%%%%%%%%%%%%%%%%%%%%%%%%%%%%%%%%%%%%%%%%%%%%%%%%%%%%%%%%%%
This section is devoted to the study of the auxiliary curve mentioned at the end of the
previous section. This curve is obtained computing a Gr\"{o}bner basis to eliminate
$x_1,x_2$ and $u$ in the system $\mathfrak{S}_2(d,k)$. Doing this elimination, one
arrives at the following definition:
\begin{defn}\label{defCurvaS}
Let $s$ be the polynomial:
\[s(\bar y,d,k)=(f_2^2+f_1^2)(y_2-k)^2-f_2^2d^2.\]
For every $(d_0,k_0)\in{\mathbb C}^2$, the {\sf auxiliary curve} ${\mathcal S(d_0,k_0)}$
to ${\mathcal C}$ is the affine plane curve defined over ${\mathbb C}$ by the polynomial
$s(\bar y,d_0,k_0)$.
\end{defn}

The following theorem relates the solutions in Theorem \ref{Theorem01} with the
intersection points of ${\mathcal C}$ and the auxiliary curve.
\begin{thm}\label{Theorem02}
Let $\Delta$ be as in Theorem \ref{Theorem01}, let $(d_0,k_0)\in\Delta$, and let
$\Gamma$ be the set of $\delta_1$ solutions of ${\mathfrak S}_2(d_0,k_0)$ appearing in
Theorem \ref{Theorem01}. Then it holds that:
\begin{enumerate}
\item[(a)] If $(p,q,u_0)\in\Gamma$, then $p\in{\mathcal C}\cap{\mathcal S(d_0,k_0)}$.
\item[(b)] If $p\in{\mathcal C}\cap{\mathcal S(d_0,k_0)}$ and $p$ is not of ramification
of ${\mathcal C}$, there exist $q\in{\mathbb C}^2$ and $u_0\in{\mathbb C}$ such that
$(p,q,u_0)\in\Gamma$.
\end{enumerate}
\end{thm}

\begin{rem}
\begin{enumerate}
\item[]
\item The solution $(p,q,u_0)$ in statement (b) of Theorem \ref{Theorem02} can be expressed as:
\[a_1=\dfrac{(-f_1(p)b_2+f_2(p)b_1+f_1(p)k_0)}{f_2(p)}\,,\quad a_2=k_0,\qquad u_0=\dfrac{1}{f_1^2(p)+f_2^2(p)},\]
where $p=(b_1,b_2)$ and $q=(a_1,a_2)$. Also note that, since $p$
is not of ramification of ${\mathcal C}$, it follows that $f_2(p)\neq 0$. Moreover, since $f_2(p)\neq 0, d_0\neq 0$ and $p\in{\mathcal S(d_0,k_0)}$, then ${f_1^2(p)+f_2^2(p)}\neq 0$.
\item  Note that ${\mathcal C}\cap{\mathcal S(d_0,k_0)}$ may contain other points besides those appearing
in the theorem. For example, every affine singularity of ${\mathcal C}$ is also a point
of ${\mathcal C}\cap{\mathcal S(d_0,k_0)}$. But the theorem shows a 1:1 correspondence
between $\Gamma$ and the points in ${\mathcal C}\cap{\mathcal S(d_0,k_0)}$ that are not
of ramification in $\mathcal C$.
\end{enumerate}
\end{rem}
{\bf Proof. }{\small

(a) We consider the polynomials
\[
\begin{cases}
\nu_1(\bar y)=-f_2^2(\bar y)\\
\nu_2(\bar y,\bar x)=f_1(\bar y)(x_2-y_2)+f_2(\bar y)(x_1-y_1)\\
\nu_3(\bar y,\bar x,k)=(f_2^2(\bar y)+f_1^2(\bar y))(2y_2-x_2-k))
\end{cases}
\]
Then it can be easily checked that
\[s(\bar y,d,k)=\nu_1(\bar y)b(\bar y,\bar x, d)+\nu_2(\bar y,\bar x) n(\bar y,\bar x)+
\nu_3(\bar y,\bar x)\ell(\bar y,k)
\]
Now, let $(p,q,u_0)\in\Gamma$. Then by Theorem \ref{Theorem01}(1b), one has that
$p\in\mathcal C$. Moreover, because of the above description of the polynomial $s$, and
taking into account
 that $(p,q,u_0)$ is a solution of $\mathfrak{S}_2(d_0,k_0)$, one has that $p\in\mathcal S(d_0,k_0)$.

(b) Let $p=(b_1,b_2)\in\mathcal C\cap\mathcal S(d_0,k_0)$ be such that $f_2(p)\neq 0$.
Then we consider
\[q=(a_1,a_2)=\left(\dfrac{-f_1(p)b_2+f_2(p)b_1+f_1(p)k_0}{f_2(p)},k_0\right)\]
and
\[u_0=\dfrac{1}{f_1^2(p)+f_2^2(p)}.\]

Note that $s(p,d_0,k_0)=(f_1^2(p)+f_2^2(p))(b_2-k_0)^2-f_2^2(p)d_0^2=0$ and
$f_2(p)d_0^2\neq 0$, and hence $f_1^2(p)+f_2^2(p)\neq 0$. Now, let us see that
$(p,q,u_0)\in\Gamma$. Substituting $(p,q,u_0)$ in $\mathfrak S_2(d_0,k_0)$ one sees that
it is a solution of the system. Moreover, $p\in\mathcal C$, it is regular and it is not
of ramification. Furthermore, because of the vanishing of $f, b, n$ and $\ell$ at
$(p,q,u_0)$, one has that $q\in\mathcal L(k_0)\cap\mathcal O_{d_0}(\mathcal C)$.
Therefore $(p,q,u_0)\in\Gamma$.}

In Theorem \ref{Theorem02} we have seen that (generically in $(d,k)$) there is a 1:1
correspondence between the $\delta_1$ points in $\Gamma$ and the points in ${\mathcal
C\cap S(d,k)}$ where $f_2$ does not vanish. The advantage of this strategy is that,
while the generic offset equation is not known, both $f$ and $s$ are known polynomials.
Therefore we can use standard techniques, such as those provided by B\'{e}zout's Theorem, to
analyze the intersection points between the two plane curves. But, for our purposes, we
have to ensure the following: first, we are going to consider all the intersection
points of ${\mathcal C}$ and ${\mathcal S}(d,k)$, so we have to treat the problem
projectively. Thus, we consider the projective closures of the curves, and we {\bf
denote} them by $\overline{\mathcal C}$ and $\overline{\mathcal S(d,k)}$, respectively.
Secondly, $\overline{\mathcal C}\cap\overline{\mathcal S(d,k)}$ may contain also points
that are not associated to points in $\Gamma$, and we need to distinguish them. This
fact motivates the following definition.

\begin{defn}\label{def01}\label{defFakePoints}
Let $\Delta$ be as in Theorem \ref{Theorem01}, and let $(d_0,k_0)\in\Delta$.
\begin{enumerate}
\item[]
\item The affine intersection points of $\overline{\mathcal C}$ and $\overline{\mathcal
S(d_0,k_0)}$ that are not of ramification of $\overline{\mathcal C}$ are called {\sf
non-fake points}.
\item The remaining intersection points of $\overline{\mathcal C}$ and
$\overline{\mathcal S(d_0,k_0)}$ are called {\sf fake points}.
\end{enumerate}
We denote by ${\mathcal F}$ the set of all fake points.
\end{defn}
\begin{rem}\label{NumberNonFakePoints}
Observe that because of Theorems \ref{Theorem01} and \ref{Theorem02}, for each
$(d_0,k_0)\in\Delta$ the number of non-fake points is precisely the partial degree
$\delta_1$.
\end{rem}

Although $\mathcal F$ seems to depend on the choice of $(d_0,k_0)\in\Delta$, in the next
proposition we show that it is in fact invariant. Nevertheless, the set of non-fake
points does depend on $(d,k)$. Since we are working projectively, we {\bf denote} by $F,
F_1, F_2$ and $S$ the homogenization w.r.t. a new variable $y_3$ of the
polynomials $f, f_1, f_2$ and $s$ respectively. We also {\bf denote} $\bar
y_H=(y_1:y_2:y_3)$. Observe that:
\[S=(F_2^2+F_1^2)(y_2-ky_3)^2-F_2^2y_3^2d^2.\]

\begin{prop}[\sf Invariance of the fake points]\label{CaracterizacionInvarianteFakes}
The set $\mathcal F$ is finite, and does not depend on $\{d,k\}$. Furthermore,
$p\in\mathcal F$ if and only if $p\in\overline{\mathcal C}$ and either $p$ is affine and
singular or $p$ is $(1:0:0)$ or $p$ is at infinity satisfying $F_1^2(p)+F_2^2(p)=0$.
\end{prop}
{\bf Proof. }{\small

Let $p=(a:b:c)\in\mathcal F$. Then there exists $(d_0,k_0)\in\Delta$ ($\Delta$ as in
Theorem \ref{Theorem01}), such that $p\in\overline{\mathcal C}\cap\overline{\mathcal
S(d_0,k_0)}$ and either $c\neq 0$ and $F_2(p)=0$ or $c=0$. If $c=0$, since
$S(p,d_0,k_0)=0$ one has that $(F_1^2(p)+F_2^2(p))b=0$, and hence either $p=(1:0:0)$ or
$p$ is at infinity and it is isotropic. On the other hand, if $c\neq 0$ and $F_2(p)=0$,
since $p\in\overline{\mathcal S(d_0,k_0)}$ one has that $F_1(p)(b-k_0c)=0$. Now, because
of the construction of $\Delta$ (see how $\Delta_4$ is defined in the proof of Theorem
\ref{Theorem01}), $b-k_0c\neq 0$. Therefore, $p$ is affine and singular.

Conversely, if $p\in\overline{\mathcal C}$ and it satisfies any of the three conditions
in the statement of the proposition, then $p\in\overline{\mathcal S(d_0,k_0)}$. Thus, by
Definition \ref{defFakePoints} the implication holds.

Finally, from the above characterization it follows that $\mathcal F$ is
finite.\hfill\qed}

\begin{rem}\label{RemarkPlaces}
Let $p=(a:b:1)$ be a non-fake point. Observe then that necessarily $b-k_0\neq 0$, for
every $(d_0,k_0)\in\Delta$ (see the proof of Proposition
\ref{CaracterizacionInvarianteFakes}).
\end{rem}

In order to apply B\'{e}zout's Theorem we need to prove that $\overline{\mathcal C}$ and
$\overline{\mathcal S(d_0,k_0)}$ do not have common components, and we have to analyze
 the multiplicity of intersection of $\overline{\mathcal C}$ and $\overline{\mathcal
S(d_0,k_0)}$ at the non-fake points. This is the content of the following proposition:
\begin{prop}[\sf B\'{e}zout's Theorem preparation]\label{Prop02}
There exists a non-empty open subset $\tilde\Delta\subset\Delta$, where $\Delta$ is as
in Theorem \ref{Theorem01}, such that for every $(d_0,k_0)\in\tilde\Delta$ the following
hold:
\begin{enumerate}
\item  $\deg({\mathcal S}(d_0,k_0))=2\deg({\mathcal C})$,
\item ${\mathcal C}$ and ${\mathcal S(d_0,k_0)}$ have no common component,
\item  if $p$ is a non-fake point,
then $\operatorname{mult}_p({\mathcal C},{\mathcal S}(d_0,k_0))=1$.
\item Let $S(\bar y_H,d,k)$ be considered as an element of $(\mathbb
C[\bar y_H])[d,k]$:
\[S(\bar y_H,d,k)=Z_{2,0}(\bar y_H)d^2+Z_{0,2}(\bar y_H)k^2+Z_{0,1}(\bar y_H)k+Z_{0,0}(\bar y_H)\]
where:
\[
\begin{cases}
Z_{2,0}(\bar y_H)=-F_{2}^{2}y_3^2 \\
Z_{0,2}(\bar y_H)=(F_2^2+F_1^2)y_3^2,\\
Z_{0,1}(\bar y_H)=-2(F_2^2+F_1^2)y_2y_3,\\
Z_{0,0}(\bar y_H)=(F_2^2+F_1^2)y_2^2,
\end{cases}
\]
and let $\mathcal J_{\alpha}$ be the curve defined by $Z_{\alpha}(\bar y_H)$. Then it
holds that: \[\bigcap_{\alpha}(\overline{{\mathcal C}} \cap \overline{{\mathcal
J_{\alpha}}})\subset\mathcal F.\]
\item $(0:0:1)\not\in\left(\overline{{\mathcal C}} \cap \overline{{\mathcal
S(d_0,k_0)}}\right)\setminus{\mathcal F}$
\end{enumerate}
\end{prop}

{\bf Proof. }{\small

(1) $S=(F_2^2+F_1^2)(y_2-ky_3)^2-F_2^2y_3^2d^2$. The form $F_2^2y_3^2d^2$ has degree
$2n$ in $\bar y$ for $d\neq 0$, and the form $(F_2^2+F_1^2)(y_2-ky_3)^2$ has degree less
or equal than $2n$ in $\bar y$. Thus $\deg_{\bar y}(S)\leq 2n$. Now the degree could
only drop if the two forms were identical, which is generically impossible, since $d$
does not appear in the first one and $k$ does not appear in the second one. Thus, our
claim holds.

(2) Let us see that for $(d_0,k_0)\in\Delta$, $\overline{\mathcal C}$ and
$\overline{\mathcal S(d_0,k_0)}$ have no common components. Assume that they do. Then,
since $F$ is irreducible, there exists $K(\bar y_H)\in\mathbb C[y_1,y_2,y_3]$ such that
\[S(\bar y_H,d_0,k_0)=K(\bar y_H)F(\bar y_H).\]
Now, we will see that then $F_2$ vanishes on almost all point of $\mathcal C$. That
implies that $\mathcal C$ is a line, which is impossible by assumption. Indeed, if there
were infinitely many points in $\overline{\mathcal C}\cap\overline{\mathcal S(d_0,k_0)}$
with $F_2\neq 0$, this would imply infinitely many affine points in $\overline{\mathcal
C}\cap\overline{\mathcal S(d_0,k_0)}$ with $f_2\neq 0$. Then Theorems \ref{Theorem01}
and \ref{Theorem02} would give an infinite number of affine intersections between the
line $x_2-k_0=0$ and the offset, which is impossible; note that if $\mathcal
O_{d_0}(\mathcal C)$ contains a line, then $\mathcal C$ is a line.

(3) Let $(d_0,k_0)\in\Delta$, and let $p=(a,b)$ be a non-fake point. By definition, we
know that $p$ is an affine regular point of ${\mathcal C}$. Therefore, there is only one
branch of ${\mathcal C}$ passing through $p$. Let $q$ be the point in ${\mathcal
O}_{d_0}({\mathcal C})\cap{\mathcal L(k_0)}$ associated with $p$ (see Theorem
\ref{Theorem01}, (1a) for the existence of $q$). Also, by Theorem \ref{Theorem01}(1a),
$\operatorname{mult}_q({\mathcal O}_{d_0}({\mathcal C}),{\mathcal L(k_0)})=1$. Thus it
is enough to prove that $\operatorname{mult}_p({\mathcal C},{\mathcal
S}(d_0,k_0))=\operatorname{mult}_q({\mathcal O}_{d_0}({\mathcal C}),{\mathcal L}(k_0))$.
The proof will proceed as follows:
\begin{enumerate}
\item First, we consider a place $P(t)=(y_1(t),y_2(t))$ of ${\mathcal C}$ centered at
$p$, and we compute $s(P(t))$. Note that the order of this formal power series is
$\operatorname{mult}_p({\mathcal C},{\mathcal S}(d_0,k_0))$.
\item Second, we use $P(t)$ to obtain a place $Q(t)$ of ${\mathcal O}_{d_0}({\mathcal C})$
centered at $q$, and we obtain $\ell(Q(t),k_0)$. Note that the order of this formal
power series is $\operatorname{mult}_q({\mathcal O}_{d_0}({\mathcal C}),{\mathcal
L(k_0)}).$
\item Finally we prove that $\operatorname{ord}\left(\ell(Q(t),k_0)\right)=\operatorname{ord}(s(P(t))).$
\end{enumerate}

Let
\[
\begin{cases}
f_1(P(t))=v_1+\alpha t+\cdots\\
f_2(P(t))=v_2+\beta t+\cdots
\end{cases}
\]
for some $v_1,v_2,\alpha,\beta\in\mathbb C$, where $f_1(p)=v_1, f_2(p)=v_2$. This means
that the tangent vector to ${\mathcal C}$ at $p$ is $(-v_2,v_1)$ and so, there exists
$\lambda$ such that the place $P(t)$ can be expressed in the form:
\[
P(t):\begin{cases}
y_1=a-\lambda v_2 t+\cdots\\
y_2=b+\lambda v_1 t+\cdots
\end{cases}
\]
The notation $T_0=\sqrt{v_1^2+v_2^2}$ and $T_1=v_1\alpha+v_2\beta$ will be used in the
rest of the proof. Note that, since $(d_0,k_0)\in\Delta$, and $p$ is non-fake, then
$v_2, T_0$ and $b-k_0$ are all not zero (see Remark \ref{RemarkPlaces}). Now,
substituting $P(t)$ into the polynomial $s(y_1,y_2,d_0,k_0)$ leads to a power series,
whose zero-order term coefficient $A_0$ must vanish (because $p\in{\mathcal
S}(d_0,k_0)$). This term is:
\[A_0=(v_1^2+v_2^2)(b-k_0)^2-d_0^2v_2^2=T_0^2(b-k_0)^2-d_0^2v_2^2\]
Therefore we get that:
\[T_0^2=\dfrac{-(d_0v_2)^2}{(b-k_0)^2}\]
The coefficient of the first-order term $A_1$ of $s(P(t))$ is:
\[A_1=2(-d_0^2v_2\beta+T_0^2(b-k_0)\lambda v_1+(v_1\alpha+v_2\beta)(b-k_0)^2).\]
Next, using $P(t)$, we generate a place $Q(t)$ of ${\mathcal O}_{d_0}({\mathcal C})$
centered at $q$. If $(y_1,y_2)$ is a regular point in ${\mathcal C}$, the associated
point $(x_1,x_2)$ in ${\mathcal O}_{d_0}({\mathcal C})$  is given by:
\[(x_1,x_2)=(y_1,y_2)\pm d_0\dfrac{(f_1,f_2)}{\sqrt{f_1^2+f_2^2}}\]
Moreover, since $v_1^2+v_2^2\neq 0$, the power series
\[f_1^2(P(t))+f_2^2(P(t))=(v_1^2+v_2^2)+2(v_1\alpha+v_2\beta)t+\cdots
\]
has order zero (is a unit), and hence
\[\dfrac{1}{\sqrt{f_1^2(P(t))+f_2^2(P(t))}}\]
can be expressed as the following formal power series.
\[\dfrac{1}{\sqrt{f_1^2(P(t))+f_2^2(P(t))}}=
\dfrac{1}{\sqrt{v_1^2+v_2^2}}- \dfrac{v_1\alpha+v_2\beta}{(v_1^2+v_2^2)^{3/2}}t+\cdots
\]
So:
\[\begin{cases}\dfrac{f_1}{\sqrt{f_1^2+f_2^2}}={\dfrac {v_{{1}}}{{T_0}}}+\left(
{\dfrac {\alpha}{{T_0}}}-{\dfrac {{T_1}\,v_{{1}}}{{T_0}^{3}}}\right)t+\cdots\\
\dfrac{f_2}{\sqrt{f_1^2+f_2^2}}= {\dfrac {v_{{2}}}{{T_0}}}+\left({\dfrac {\beta}{
{T_0}}} -{\dfrac {{T_1}v_{{2}}}{{T_0}^{3}}}\right)t+\cdots\end{cases}
\]
Therefore $Q(t)$ is one of the two places:
\[Q(t)=(x_1(t),x_2(t))=P(t)\pm d_0\dfrac{(f_1(P(t)),f_2(P(t)))}{\sqrt{f_1^2(P(t))+f_2^2(P(t))}},\]
and so:
$$
\begin{cases}
x_1(t)=\left(a\pm{\dfrac {d_0v_{{1}}}{{T_0}}} \right)+\left(-\lambda\,v_{{2}}\pm{\dfrac
{d_0
\alpha}{{T_0}}}\mp{\dfrac{d_0\,T_1\,v_1}{{T_0}^{3}}} \right)t +\cdots\\[3mm]
x_2(t)=\left(b\pm{\dfrac {d_0v_{{2}}}{{T_0}}} \right)+\left(\lambda\,v_{{1}}\pm{\dfrac
{d_0 \beta}{{T_0}}}\mp{\dfrac{d_0\,T_1\,v_2}{{T_0}^{3}}} \right)t +\cdots
\end{cases}
$$
Substituting $Q(t)$ in the line ${\mathcal L}(k_0)$ one has:
\[
x_2(t)-k_0=\left(b\pm{\dfrac {d_0v_{{2}}}{{T_0}}}-k_0
\right)+\left(\lambda\,v_{{1}}\pm{\dfrac {d_0
\beta}{{T_0}}}\mp{\dfrac{d_0\,T_1\,v_2}{{T_0}^{3}}} \right)t +\cdots=B_0+B_1t +\cdots
\]
Now, since $\operatorname{mult}_q({\mathcal O}_{d_0}({\mathcal C}),{\mathcal
L}(k_0))=1$, one has that
\[B_0=\left(b\pm{\dfrac {d_0v_{{2}}}{{T_0}}}-k_0 \right)=0\mbox{, and }B_1=\left(\lambda\,v_{{1}}\pm{\dfrac
{d_0 \beta}{{T_0}}}\mp{\dfrac{d_0\,T_1\,v_2}{{T_0}^{3}}} \right)\neq 0
\]
Therefore
\[\pm T_0=-\dfrac{d_0v_2}{b-k_0}\]
Substituting the above equality in $B_1$ one gets
\[B_1=\dfrac{1}{T_0^3}\left(\mp\lambda\,v_{{1}}\left(\dfrac{d_0v_2}{b-k_0}\right)^3\pm d_0\beta \left(\dfrac{d_0v_2}{b-k_0}\right)^2\mp d_0T_1v_2
\right)=\]
\[
\dfrac{\mp d_0v_2}{T_0^3(b-k_0)^3}\left(d_0^2v_2^2\lambda\,v_{{1}}- d_0^2\beta
v_2(b-k_0)+ T_1(b-k_0)^3 \right)
\]
Note that this result does not depend on the previous choice of sign. And using the same
equality in $A_1$ gives:
\[A_1=2(-d_0^2\beta v_2+\left(\dfrac{d_0v_2}{b-k_0}\right)^2(b-k_0)\lambda v_1+T_1(b-k_0)^2)=\]
\[\dfrac{2}{b-k_0}\left(-d_0^2\beta v_2(b-k_0)+d_0^2v_2^2\lambda v_1+T_1(b-k_0)^3\right)\]
We observe that the term in parenthesis in $A_1$ and $B_1$ coincides. Since $B_1\neq 0$,
one has that $A_1\neq 0$ and $\operatorname{mult}_p({\mathcal C},{\mathcal
S}(d_0,k_0))=1$.

(4) Since we have assumed that $f$ does not divide to $f_{1}^{2}+f_{2}^{2}$ (in
particular $f_1^2+f_2^2\neq 0$), and that $\mathcal C$ is not a line (in particular
$f_2\neq 0$), all $\mathcal J_\alpha$ are algebraic curves. Now
\[ \bigcap_{\alpha} [\overline{{\mathcal C}}\cap \overline{{\mathcal J}_{\alpha}}]\subset
\overline{{\mathcal C}}\cap \overline{{\mathcal J}_{(2,0)}} \] and by Proposition
\ref{CaracterizacionInvarianteFakes}, $\overline{{\mathcal C}}\cap \overline{{\mathcal
J}_{(2,0)}}\subset {\mathcal F}$.

(5) Let $p=(0:0:1)$ and $A(d,k)=S(p,d,k)$. If either $p\in\mathcal F$ or
$p\not\in{\mathcal C}$, then no further restriction on $\Delta$ is required. Now, let
$p\in\mathcal C$ and $p\not\in\mathcal F$. Then by Proposition
\ref{CaracterizacionInvarianteFakes}, $P$ is not a singularity of $\mathcal C$. Now, if
$F_2(p)\neq 0$, then $A$ is not constant. Moreover, if $F_2(p)=0$, then $F_1(p)\neq 0$
and $A$ is not constant either. Let $\Psi$ be the curve in $\mathbb C^2$ defined by $A$.
Then in $\Delta\setminus\Psi$ statement (5) holds. Indeed, if
$p\in\left(\overline{{\mathcal C}} \cap \overline{{\mathcal
S(d_0,k_0)}}\right)\setminus{\mathcal F}$, then $p\in\mathcal C$, $A(d_0,k_0)=0$ and
$p\not\in\mathcal F$. Thus $(d_0,k_0)\in\Psi$.}

%%%%%%%%%%%%%%%%%%%%%%%%%%%%%%%%%%%%%%%%%%%%%%%%%%%%%%%%%%%%%%%%%%%%%%%%%%%%%%%%%%%%%%%%%
%%%%%%%%%%%%%%%%%%%%%%%%%%%%%%%%%%%%%%%%%%%%%%%%%%%%%%%%%%%%%%%%%%%%%%%%%%%%%%%%%%%%%%%%%
\section{Cornerstone Theorem}\label{secCornerstoneTheorem}
%%%%%%%%%%%%%%%%%%%%%%%%%%%%%%%%%%%%%%%%%%%%%%%%%%%%%%%%%%%%%%%%%%%%%%%%%%%%%%%%%%%%%%%%%
%%%%%%%%%%%%%%%%%%%%%%%%%%%%%%%%%%%%%%%%%%%%%%%%%%%%%%%%%%%%%%%%%%%%%%%%%%%%%%%%%%%%%%%%%

Later, in Section \ref{SecStrategyDistance}, when analyzing the
problem of the degree in $d$ of the generic offset, we will find
another situation which involves the intersection of $\mathcal C$
with an auxiliary curve that plays the role that $S$ plays here, and
a concept of fake and non-fake intersection points with properties
analogous to those described in the previous results. The next
result shows how those properties of an auxiliary curve can be used
to establish a degree formula. We will give a general formulation in
order to apply this same result to both situations. In the statement
of the next theorem we use the following {\bf terminology}: let
$\bar u=(u_1,u_2)$. Then, if $h\in{\mathbb C}[y_1,y_2,y_3,\bar u]$,
we denote by $\operatorname{PP}_{\bar u}(h)$ the primitive part of
$h$ w.r.t. $\bar u$, and by $\operatorname{Res}_{y_3}(h_1,h_2)$ the
resultant of $h_1,h_2\in{\mathbb C}[y_1,y_2,y_3,\bar u]$ w.r.t.
$y_3$. Recall that $\bar y_H=(y_1:y_2:y_3)$.

\begin{thm}[\sf Cornerstone Theorem]\label{Theorem07} Let $\mathcal D$
be an irreducible affine plane curve, not being a line, and let $Z(\bar y_H,\bar
u)\in\mathbb C[\bar y_H,\bar u]$ be homogeneous in $\bar y_H$ and depending on $y_3$.
Let us suppose that there exists an open set $\Xi\subset\mathbb C^2$ such that, for
$\bar \omega\in\Xi$ the following hold:
\begin{enumerate}
\item $\deg_{\bar y_H}(Z(\bar y_H, \bar \omega))=\deg_{\bar y_H}(Z(\bar y_H, \bar u))$.
Let $\mathcal Z(\bar \omega)$ be the plane curve defined  by $Z(\bar y_H, \bar \omega)$
(note that $Z(\bar y_H, \bar \omega)$ is non-constant).
\item $\mathcal Z(\bar \omega)$ and $\mathcal D$ do not have common
components.
\item Let
\[ {\mathcal G}=\bigcap_{\bar u\in \Xi} [\overline{{\mathcal Z}(\bar
u)} \cap \overline{{\mathcal D}}].
\] Then, for every $p\in [\overline{{\mathcal Z}(\bar\omega)}\cap \overline{{\mathcal D}}]\setminus {\mathcal G}$, $\operatorname{mult}_p(\overline{
{\mathcal D}},\overline{{\mathcal Z(\bar \omega)}})=1$.
\item Let $Z(\bar y_H,\bar u)$ be considered as an element of $(\mathbb
C[\bar y_H])[\bar u]$, so that one has:
\[Z(\bar y_H,\bar u)=\sum_{\alpha}Z_{\alpha}(\bar y_H)\bar u^{\alpha}\]
for some $Z_{\alpha}(\bar y_H)\in\mathbb C[\bar y_H]$. If
$Z_{\alpha}(\bar y_H)$ is not constant, let $\mathcal J_{\alpha}$ be
the curve it defines. Then it holds that:
\[\bigcap_{\alpha}(\overline{{\mathcal D}}\cap \overline{\mathcal J_{\alpha}})\subset\mathcal G.\]
\item $(0:0:1)\not\in\left(\overline{{\mathcal
Z(\bar\omega)}} \cap \overline{{\mathcal D}}\right)\setminus{\mathcal G}$
\end{enumerate}
Then, there exists a non-empty open subset $\Xi^\star\subset \Xi$ such that for
$\bar\omega\in \Xi^\star$:
\[\fbox{$\operatorname{Card}([\overline{{\mathcal Z}(\bar\omega)}\cap \overline{{\mathcal D}}]\setminus {\mathcal G})=
\deg_{\{y_1,y_2\}}\left(\operatorname{PP}_{\bar u}\left( \operatorname{Res}_{y_3}(G(\bar
y_H),Z(\bar y_H,\bar u))\right)\right)$},\] where $G$ is the form defining the
projective closure $\overline{\mathcal D}$ of $\mathcal D$.
\end{thm}

{\bf Proof. }{\small

We denote by $R(y_1,y_2,\bar u)=\operatorname{Res}_{y_3}(G,Z)$; observe that, since $G$
is irreducible and $\mathcal D$ is not a line, $G$ depends on $y_3$, moreover $Z$
depends also on $y_3$ by hypothesis. Let $R(y_1,y_2,\bar u)$ factor as
$$R(y_1,y_2,\bar u)=M(y_1,y_2)N(y_1,y_2,\bar u)$$ where $M$ and $N$
are the content and primitive part of $R$ w.r.t. $\bar u$, respectively. Then $M$ and
$N$ are homogeneous polynomials in $y_1,y_2$, and $M\in {\mathbb C}[y_1,y_2], N\in
{\mathbb C}[\bar u][y_1,y_2]$. This implies that $M$ factors over $\mathbb C$ in linear
factors, namely:
$$M=\prod_{i=1}^r(\beta_iy_1-\alpha_iy_2)$$

We observe that the leading  coefficient $L$ of $Z$ w.r.t. $y_3$ is a non-zero
polynomial in ${\mathbb C}[\bar u][y_1,y_2]$. If $L$ does not depend on $\bar u$ or any
coefficient of $L$ w.r.t. $\{y_1,y_2\}$ is a non-zero constant we take $\Psi=\emptyset$,
otherwise we take $\Psi$ as  the intersection of all curves  in ${\mathbb C}^2$ defined
by each non-constant coefficient of $L$ w.r.t. $\{y_1,y_2\}$. Let $\Xi_1=\Xi \setminus
\Psi$. Since $G$ does not depend on $\bar u$, for every $\bar \omega\in \Xi_1$, both
leading coefficients of $G$ and $Z(\bar{y}_H,\bar \omega)$ w.r.t. $y_3$ do not vanish.
In particular, this implies that the resultant specializes properly; i.e. if
$Z_0(\bar{y}_{H})=Z(\bar{y}_H,\bar \omega)$ and
$R_0(y_1,y_2)=\operatorname{Res}_{y_3}(G,Z_0)$, then for $\bar\omega\in \Xi_1$
$$R_0=M(y_1,y_2)N(y_1,y_2,\bar \omega).$$
By Lemma 18 in \cite{SSS05}, and because of $\Xi_1$ and hypothesis (1), we observe that
$R$ and $R_0$ have the same degree. Hence the degree of $N(y_1,y_2,\bar u)$ and
$N_0=N(y_1,y_2,\bar\omega)$ is also the same. Moreover, since $N_0$ is a homogeneous
polynomial, it can be factored as
$$N_0=\prod_{j=1}^s(\beta'_jy_1-\alpha'_jy_2).$$
Thus
$$R_0=M\cdot N_0=\prod_{i=1}^r(\beta_iy_1-\alpha_iy_2)\prod_{j=1}^s(\beta'_jy_1-\alpha'_jy_2)$$
In this situation, for $\bar \omega\in \Xi$ let ${\mathcal B}_{\bar
\omega}=[\overline{{\mathcal Z}(\bar \omega)}\cap \overline{\mathcal D}]\setminus
\mathcal G$. Then, since $\deg(N)=\deg(N_0)$, the proof ends if we find a non-empty open
subset $\Xi^\star\subset \Xi$ such that $\operatorname{Card}({\mathcal B}_{\bar
\omega})=\deg(N_0)$ for $\bar \omega \in \Xi^\star$.

We start the construction of $\Xi^\star$. First, we prove that there exists a non-empty
open subset $\Xi_2\subset\Xi_1$ such that, if $\bar\omega\in\Xi_2$, then
$\gcd(N_0,M)=1$. Indeed, first we observe that $\gcd(N,M)=1$, since otherwise $N$ would
have a factor depending on $\{y_1,y_2\}$, and $N(y_1,y_2,\bar u)$ is primitive w.r.t.
$\bar u$. Now, for each factor $(\beta_i y_1-\alpha_i y_2)$ of $M$, we consider the
polynomial $N(\alpha_i,\beta_i,\bar u)$. This polynomial is not identically zero because
$\gcd(N,M)=1$. Then $\Xi_2=\Xi_1\setminus (\Gamma_1 \cup \cdots \cup \Gamma_r)$, where
$\Gamma_i$ is the curve in ${\mathbb C}^2$ defined by $N(\alpha_i,\beta_i,\bar u)$.

Now, we prove the existence of a non-empty open subset $\Xi_3\subset \Xi_2$ such that
for $\bar \omega \in \Xi_3$
 the projective lines $\overline{{\mathcal L}_i}$, defined by the equations $\beta_i
y_1-\alpha_i y_2=0$, do not contain points of ${\mathcal B}_{\bar\omega}$; recall that
$\beta_1y_1-\alpha_i y_2$ is a factor of $M$. For this purpose, observe that
$\overline{{\mathcal L}_i}$ meets $\overline{\mathcal D}$ in a finite number of points;
recall that by assumption $\mathcal D$ is irreducible and it is not a line. Let $
[\overline{{\mathcal D}}\cap \overline{{\mathcal L}_i}]\setminus {\mathcal
G}=\{P_{1}^{i},\ldots,P_{k_i}^{i}\}.$
 Now, consider the polynomials
$Z(P_{j_i}^{i},\bar u)$. These polynomials are not identically zero, because otherwise
it would imply that all coefficients of $Z(\bar{y}_{H},\bar u)$ w.r.t. $\bar u$ vanish
at $P_{j_i}$, and by hypothesis (5), that
\[ P_{i_j}^{i}\in \bigcap_{\alpha} (\overline{\mathcal D}\cap \overline{{\mathcal
J}_{\alpha}}) \subset {\mathcal G}, \] which is impossible. Then, if $\Psi_{j_i}^{i}$ is
the curve in ${\mathbb C}^2$ defined by $Z(P_{j_i}^{i},\bar u)$, Let \[ \Xi_3=\Xi_2
\setminus [\bigcup_{i=1}^{r} \bigcup_{j=1}^{k_i} \Psi_{j_i}^{i}].
\]
Let us see that $\Xi_3$ satisfies the requirements. Let $\bar \omega \in \Xi_3$, and
assume that there exists  $P\in [\overline{{\mathcal L}_i}\cap \overline{{\mathcal
Z}(\bar \omega)}\cap \overline{{\mathcal D}}] \setminus {\mathcal G}$. Then, $P\in
[\overline{{\mathcal L}_i} \cap \overline{{\mathcal D}}]\setminus {\mathcal G}$.
Therefore there exists $j_i$ such that $P=P_{j_i}^{i}$, and because $P\in
\overline{{\mathcal Z}(\bar \omega)}$ one has that $Z(P_{j_i}^{i},\bar \omega)=0$, which
is a contradiction since $\bar\omega \not\in \Psi_{j_i}^{i}$.

Finally the last open subset is constructed. Let $W(y_1,y_2)$  be the leading
coefficient of $G(\bar{y}_{H})$ w.r.t. $y_3$. Note that $W\in {\mathbb C}[y_1,y_2]$ is
homogeneous. Then, we choose a non-empty Zariski open subset $\Xi_4\subset\Xi_3$ such
that for every $\bar \omega \in \Xi_4$ it holds that $\gcd(N_0,W)=1$. For this purpose,
let $W$ factor as
\[ W=\prod_{k=1}^{m} (\sigma_i y_1-\nu_i y_2). \]
We consider the polynomials $N(\nu_i,\sigma_i,\bar{u})$. These polynomials are not
identically zero, because otherwise it would imply (note that $N$ is homogeneous in
$y_1,y_2$) that $N$ has a factor, namely $(\sigma_i y_1-\nu_i y_2)$, and $N$ is
primitive w.r.t. $\bar u$. Then, we consider \[ \Xi_4=\Xi_3 \setminus (\Psi_1\cup \cdots
\cup \Psi_n),\] where $\Psi_i$ is the curve in ${\mathbb C}^2$ defined by
$N(\nu_i,\sigma_i,\bar{u})$. Let us see that $\Xi_4$ satisfies the requirements. Let us
assume that $\bar \omega \in \Xi_4$ and that there exists a factor
$\Lambda=\beta'_jy_1-\alpha'_jy_2$ of $N_0=N(y_1,y_2,\bar\omega)$ such that
$\gcd(\Lambda,W)\neq 0$. Then, there exists $i\in \{1,\ldots,m\}$ such that
$\Lambda=\sigma_i y_1-\nu_i y_2$. Thus $N(\nu_i,\sigma_i,\bar \omega)=1$. That is,
$\omega\in \Psi_i$ which is a contradiction.

Now, we take $\Xi^\star=\Xi_4$, and we prove that for every $\bar\omega\in \Xi^\star$,
$\operatorname{Card}({\mathcal B}_{\bar \omega})=\deg(N_0)$:
\begin{enumerate}
\item[(a)] Let us see that if $P=(a:b:c)\in\mathcal G\setminus\{(0:0:1)\}$ then
$(by_1-ay_2)$ divides $M$. Indeed: $P\in \overline{{\mathcal Z}(\bar\omega)}\cap
\overline{\mathcal D}$ for every $\bar\omega\in \Xi^\star$. Thus, $R_0(a,b,\bar
\omega)=0$ for every $\bar\omega\in \Xi^\star$. Since the resultant specializes properly
in $\Xi^*$, because of $\Xi_1$, then $R(a,b,\bar u)=M(a,b)N(a,b,\bar u)$ vanishes on
$\Xi^\star$. Moreover, $N(a,b,\bar u)$ cannot vanish on $\Xi^\star$, since otherwise it
would imply  that $(by_1-ay_2)$ divides $N$, and $N$ is primitive w.r.t. $\bar u$. Thus,
$M(a,b)=0$.
\item[(b)] Let us see that every linear factor of $N_0$ (for every $\bar\omega\in
\Xi^\star$) generates a point in ${\mathcal B}_{\omega}$. Indeed: let $(by_1-ay_2)$
divide $N_0$ then, because of $\Xi_4$, there exists $c$ such that $(a:b:c)\in
\overline{{\mathcal Z}(\bar\omega)}\cap \overline{\mathcal D}$. Note that
$(a:b:c)\neq(0:0:1)$. Now, taking into account (a), and because of $\Xi_2$, one has that
$(a:b:c)\in {\mathcal B}_{\bar\omega}$.
\item[(c)] Let us see that every point in ${\mathcal
B}_{\bar\omega}$ (for every $\bar\omega\in \Xi^\star$) generates a factor in $N_0$.
Indeed, let $P=(a:b:c)\in {\mathcal B}_{\bar\omega}$, then by hypothesis (5)
$A=(by_1-ay_2)\neq 0$. Thus, $A$ divides $R_0$, and because of $\Xi_3$, $A$ does not
divide $M$. Therefore, $A$ divides $N_0$. \item[(d)] Now the result follows from Lemma
19 in \cite{SSS05}, from (b), (c), from hypothesis (4), and because $\gcd(M,N_0)=1$ in
$\Xi^\star$.
\end{enumerate}}

%%%%%%%%%%%%%%%%%%%%%%%%%%%%%%%%%%%%%%%%%%%%%%%%%%%%%%%%%%%%%%%%%%%%%%%%%%%%%%%%%%%%%%%%%
%%%%%%%%%%%%%%%%%%%%%%%%%%%%%%%%%%%%%%%%%%%%%%%%%%%%%%%%%%%%%%%%%%%%%%%%%%%%%%%%%%%%%%%%%
\section{Partial degree formulae for the implicit case}\label{secPartialDegreeFormulae}
%%%%%%%%%%%%%%%%%%%%%%%%%%%%%%%%%%%%%%%%%%%%%%%%%%%%%%%%%%%%%%%%%%%%%%%%%%%%%%%%%%%%%%%%%
%%%%%%%%%%%%%%%%%%%%%%%%%%%%%%%%%%%%%%%%%%%%%%%%%%%%%%%%%%%%%%%%%%%%%%%%%%%%%%%%%%%%%%%%%

Using the previous results, we derive the first two partial degree formulae for offset
curves. For the first formula we observe that, by Proposition \ref{Prop02}, and by
B\'{e}zout's Theorem, we know that for $(d_0,k_0)\in\Delta$ (with $\Delta$ as in Theorem
\ref{Theorem01})
\[\deg({\mathcal
C})\deg({\mathcal S(d_0,k_0)})= \sum_{p\in{\overline{\mathcal C}\cap\overline{\mathcal
S(d_0,k_0)}}}\operatorname{mult}_p(\overline{\mathcal C},\overline{\mathcal
S(d_0,k_0)})=\]
\[\sum_{p\in{{\mathcal F}}}\operatorname{mult}_p(\overline{\mathcal
C},\overline{\mathcal S(d_0,k_0)})+
 \sum_{p\in{\left({\overline{\mathcal C}\cap\overline{\mathcal
S(d_0,k_0)}}\right)\setminus{\mathcal F}}}\operatorname{mult}_p(\overline{\mathcal
C},\overline{\mathcal S(d_0,k_0)})
\]
Moreover, since there are $\delta_1$ non-fake points (see Remark
\ref{NumberNonFakePoints}), and for each of them the multiplicity of intersection is
one, the following formula holds.
\begin{thm}[\sf First partial degree formula]\label{Theorem03}
Let $\tilde\Delta$ be as in Theorem \ref{Prop02}. For every $(d_0,k_0)\in\tilde\Delta$,
it holds that:
$$
\fbox{
$
\delta_1=\deg_{x_1}({\mathcal O}_{d_0}({\mathcal C}))=2\left(\deg({\mathcal C})\right)^2-
\sum_{p\in{\mathcal F}}\operatorname{mult}_p(\overline{\mathcal C},\overline{\mathcal
S(d_0,k_0)})
$}$$
\end{thm}
The above formula is, although algorithmically applicable, mainly of theoretically
interest, and probably not so useful in practice, because it requires an explicit
description of the inequalities defining the open set $\Delta$.

In order to overcome this difficulty, we present a second formula that uses a univariate
resultant and gcds computations. This formula is a direct consequence of Theorem
\ref{Theorem07}. Recall that $\operatorname{PP}_{\bar u}(h)$ is the primitive part of
$h$ w.r.t. $\bar u$, and $\operatorname{Res}_{y_3}(h_1,h_2)$ is the resultant of
$h_1,h_2\in{\mathbb C}[y_1,y_2,y_3,\bar u]$ w.r.t. $y_3$. Recall also that $\bar
y_H=(y_1:y_2:y_3)$. The second partial degree formula is then the following:

\begin{thm}[\sf Second partial degree formula]\label{Theorem04}
\[
\fbox{$\delta_1=\deg_{x_1}({\mathcal O}_d({\mathcal C}))=
\deg_{\{y_1,y_2\}}\left(\operatorname{PP}_{\{d,k\}}\left(\operatorname{Res}_{y_3}(F(\bar
y_H),S(\bar y_H,d,k))\right)\right)$}\] We recall that $F$ is the homogeneous implicit
equation of $\overline{\mathcal C}$, and $S$ is the homogenization of the polynomial
introduced in Definition \ref{defCurvaS}.
\end{thm}
{\bf Proof of Theorem . }{\small In order to prove the theorem, we apply Theorem
\ref{Theorem07}. Let $\mathcal D=\mathcal C$, $Z(\bar y_H,\bar u)=S(\bar y_H,d,k)$,
where $\bar u=(d,k)$, and $\Xi=\tilde\Delta$, where $\tilde\Delta$ is as in Proposition
\ref{Prop02}. We check that all the hypothesis are satisfied:
\begin{itemize}
\item $\mathcal C$ is irreducible and it is not a line by assumption.
\item $S$ can be written as
\[S=\left((F_1^2+F_2^2)k^2-F_2^2d^2\right)y_3^2-2k(F_1^2+F_2^2)y_3+(F_1^2+F_2^2)y_2^2\]
Thus, since $F_1^2+F_2^2$ and $F_2^2$ are not identically zero, $S$ depends on $y_3$.
\item (1) and (2) in Theorem \ref{Theorem07} follow from (1) and (2)  in
Proposition \ref{Prop02}.
\item Let us see that
\[ {\mathcal F}=\bigcap_{(d,k)\in\tilde\Delta} \left[\overline{{\mathcal S}(d,k)} \cap \overline{{\mathcal C}}\right].\]
Indeed, the left-right inclusion follows from Definition \ref{defFakePoints} and
Proposition \ref{CaracterizacionInvarianteFakes}. Now, let
$p\in\bigcap_{(d,k)\in\tilde\Delta} [\overline{{\mathcal S}(d,k)} \cap
\overline{{\mathcal C}}]$. Then $p\in\mathcal C$ and $S(p,d,k)$ vanishes on $\Delta$.
Thus $S(p,d,k)$ is identically zero. So, $p\in\bigcap_{\alpha}\left( \overline{{\mathcal
C}}\cap\overline{{\mathcal J}_{\alpha}}\right)$, where ${\mathcal J}_{\alpha}$ is as in
Proposition \ref{Prop02}. Then, by Proposition \ref{Prop02}(4), one has that
$p\in\mathcal F$.
\item In this situation, hypothesis (3), (4) and (5) in Theorem \ref{Theorem07} follows from
Proposition \ref{Prop02}(3), (4) and (5), respectively.
\end{itemize}
Then, Theorem \ref{Theorem07} implies that there exists a non-empty open
$\Delta^*\subset\tilde\Delta$ such that for $(d_0,k_0)\in\Delta^*$
\[\operatorname{Card}([\overline{{\mathcal S}(d_0,k_0)}\cap \overline{{\mathcal C}}]\setminus {\mathcal F})=
\deg_{\{y_1,y_2\}}\left(\operatorname{PP}_{\{d,k\}}\left(\operatorname{Res}_{y_3}(F(\bar
y_H),S(\bar y_H,d,k))\right)\right)
\]
Now the theorem follows from Remark \ref{NumberNonFakePoints} and Proposition
\ref{CaracterizacionInvarianteFakes}. \hfill\qed}

%%%%%%%%%%%%%%%%%%%%%%%%%%%%%%%%%%%%%%%%%%%%%%%%%%%%%%%%%%%%%%%%%%%%%%%%%%%%%%%%%%%%%%%%%
%%%%%%%%%%%%%%%%%%%%%%%%%%%%%%%%%%%%%%%%%%%%%%%%%%%%%%%%%%%%%%%%%%%%%%%%%%%%%%%%%%%%%%%%%
\section{Partial degree formulae for the parametric case}\label{secParametricCase}
%%%%%%%%%%%%%%%%%%%%%%%%%%%%%%%%%%%%%%%%%%%%%%%%%%%%%%%%%%%%%%%%%%%%%%%%%%%%%%%%%%%%%%%%%
%%%%%%%%%%%%%%%%%%%%%%%%%%%%%%%%%%%%%%%%%%%%%%%%%%%%%%%%%%%%%%%%%%%%%%%%%%%%%%%%%%%%%%%%%
The formulae derived in the previous sections are valid for the implicit representation
of any irreducible algebraic plane curve. In this section, we will present a simpler
formula, adapted to the case of rational algebraic plane curves given parametrically.
This formula only requires the computation of the degree of three univariate gcds,
directly related to the parametrization.

Let
\[{\mathcal P}(t)=\left(\dfrac{X(t)}{W(t)},\dfrac{Y(t)}{W(t)}\right)\]
be a proper rational parameterization of a plane curve ${\mathcal C}$, where
\[\gcd(X,Y,W)=1.\]
As a normal vector associated to $\mathcal P(t)$ we consider $(N_1(t),N_2(t))$, where
\[
\begin{cases}
N_1(t)=-(W(t)Y'(t)-W'(t)Y(t))\\
N_2(t)=W(t)X'(t)-W'(t)X(t)
\end{cases}\]

Now, substituting in system ${\mathfrak S}_2(d,k)$ the variables $\bar y$ by the
parametrization and the partial derivatives $f_i$ by the normal vector components $N_i$,
and clearing up denominators, one may apply a similar strategy to derive the partial
degree formulae. More precisely, the auxiliary curve $\mathcal S$ is replaced here by a
univariate polynomial $\hat S(t)$ that takes values in the parameter space, namely
\[\hat S(t)=(N_1^2+N_2^2)(Wk-Y)^2-d^2W^2N_2^2.\]
A similar argument to the implicit case, based on the genericity of $k$ and $d$, shows
that the partial offset degree is the degree of the primitive part of $\hat S$ w.r.t.
$\{d,k\}$. That is:

 \fbox{$\delta_1=\deg_{x_1}({\mathcal O}_d({\mathcal C}))=
 \deg_{t}\left(\operatorname{PP}_{\{k,d\}}\left((N_1^2+N_2^2)(Wk-Y)^2-d^2W^2N_2^2\right)\right)$}

Collecting the coefficients of $\hat S$ w.r.t. $\{d,k\}$ one deduces that the content is
given by the following gcd:
\[\Theta(t)=\gcd\left(W^2\gcd(N_1,N_2)^2,(N_1^2+N_2^2)Y\gcd(W,Y)\right)\]
Since the degree of $\hat S$ equals
$2(\max(\deg(Y),\deg(W))+\max(\deg(N_1),\deg(N_2)))$, one gets the following second
formula:

\fbox{$\delta_1=2(\max(\deg(Y),\deg(W))+\max(\deg(N_1),\deg(N_2)))-\deg_{t}(\Theta(t))$}

\section{Strategy description for the distance degree formula}\label{SecStrategyDistance}

Since the generic offset equation $g$ also depends on $d$, it is natural to complete
this degree analysis by studying the degree of $g$ in $d$. We {\bf denote} it  by
$\delta_d$. We begin recalling that, for all but a finite (possibly empty) set of values
of $d$, the generic offset equation specializes properly (see Theorem \ref{Prop01}).
This implies that there are infinitely many values $d_0$ such that $g(\bar x,d_0)=0$ is
the equation of $\mathcal O_{d_0}(\mathcal C)$ and, simultaneously, $g(\bar x,-d_0)=0$
is the equation of $\mathcal O_{-d_0}(\mathcal C)$. But, because of the symmetry in the
construction, the offsets $\mathcal O_{d_0}(\mathcal C)$ and $\mathcal O_{-d_0}(\mathcal
C) $ are exactly the same. Thus, it follows that for infinitely many values of $d_0$ it
holds that up to multiplication by a non-zero constant:
\[g(\bar x,d_0)=g(\bar x,-d_0).\]
Hence, we have proved the following proposition:
\begin{prop}\label{Prop03}
The generic offset equation belongs to $\mathbb C[\bar x][d^2]$. That is, it only contains even powers of $d$. In particular, $\delta_d$
is even.
\end{prop}
\begin{rem}\label{RemarkMu}
In the sequel we {\bf denote} $\delta_d=2\mu$, where $\mu\in\mathbb N$.
\end{rem}

Now, the strategy is slightly different to the one described in Section
\ref{SecStrategyPartial}, but follows a similar structure. Essentially, it consists in
the following steps:
\begin{enumerate}
\item First, we recall that
\[n(\bar y, \bar x)=-f_2(\bar y)(x_1-y_1)+f_1(\bar y)(x_2-y_2),\]
and let
\[N(\bar y_H, \bar x)=-F_2(\bar y_H)(x_1y_3-y_1)+F_1(\bar y_H)(x_2y_3-y_2)\]
be the homogenization of $n(\bar y, \bar x)$ w.r.t. $\bar y$. For
$\bar\tau=(\tau_1,\tau_2)\in\mathbb C^2$ we {\bf denote} by ${\mathcal N}(\bar\tau)$ the
curve defined by $N(\bar y_H,\tau)$ (observe that there exists an open subset of values
of $\bar\tau$ such that ${\mathcal N}(\bar\tau)$ is indeed a curve). Let
$\overline{{\mathcal N}(\bar\tau)}$ denote the projective closure of ${\mathcal
N}(\bar\tau)$. This curve $\overline{\mathcal N}(\bar\tau)$ will play the role of the
curve $\overline{\mathcal Z}(\bar y_H,\bar u)$ used in the Cornerstone Theorem
\ref{Theorem07}.
\item Secondly, we consider the system
\[
\left.\begin{array}{lr}
F(\bar y_H)=0\\
N(\bar y_H,\bar x)=0
\end{array}\right\}
\equiv\mathfrak{S}_3(\bar\tau)
\]
and we analyze its solutions; this is done in  Theorem \ref{Theorem06}.
\item Based on this analysis, the notion of $d$-fake and non $d$-fake points are
introduced.
\item Next, the invariance of the set of $d$-fake points is established in Proposition \ref{Lemma02a}.
\item In order to apply B\'{e}zout's Theorem, we state Proposition \ref{Prop04}, which is similar to Proposition
\ref{Prop02}.
\item Finally, we apply the cornerstone Theorem.
\end{enumerate}

The second step is the content of the following theorem (compare to Theorem
\ref{Theorem01} and Theorem \ref{Theorem02}).
\begin{thm}\label{Theorem06}
There exists a non-empty Zariski open subset $U$ of $\mathbb C^2$, such that for
$\bar\tau=(\tau_1,\tau_2)\in U$
\begin{enumerate}
\item Let $\hat p$ be an affine regular point of $\mathcal C$. If $\hat p$ is the
origin or it is isotropic in $\mathcal C$, then it is not a solution of
$\mathfrak{S}_3(\bar\tau)$.
\item There exist exactly $\mu$ solutions (see Remark \ref{RemarkMu})
$\hat\Gamma(\bar\tau)=\{\hat p_i\}_{i=1,\ldots,\mu}$ of $\mathfrak{S}_3(\bar\tau)$
satisfying that $\hat p_1,\ldots,\hat p_{\mu}$ are different affine and non-isotropic
points of ${\mathcal C}$.
\item For every $\hat p_i=(a_i:b_i:1)\in\hat\Gamma(\bar\tau)$, let
\[d_i^2=(a_i-\tau_1)^2+(b_i-\tau_2)^2.\]
Then $d_1,\ldots,d_{\mu}$ are all different and non-zero.
\item For every $\hat p_i\in\hat\Gamma(\bar\tau)$, and its corresponding $d_i$ introduced in (3), it holds that $\bar\tau\in{\mathcal O}_{\pm d_i}({\mathcal
C})$, and it is the point on the offset generated by $\hat p_i  $.
\end{enumerate}
\end{thm}

{\bf Proof. }{\small  The open set $U$ is constructed in a finite number of steps, as
follows:
\begin{enumerate}
\item[(i)] Since $g$ is primitive w.r.t $d$, $g(\bar x,0)$ cannot be identically zero. Let
$\Psi_0$ be the zero set in ${\mathbb C}^2$ of $g(\bar x,0)$. And let  $U_0={\mathbb
C}^2\setminus(\mathcal C\cup\Psi_0)$.
\item[(ii)] The next open subset ensures that $\deg_{y_3}(N)$ stays invariant when specializing $\bar x$.
First, observe that none of $F_1, F_2$ cannot be identically zero because $\mathcal C$
is irreducible and it is not a line. Now, we introduce the polynomial $Z_i(y_1,y_2)$ as
the leading coefficient of $F_i$ w.r.t. $y_3$ if $F_i$ depends on $y_3$, and otherwise
$Z_i=F_i$. Let $A(\bar x,y_1,y_2)$ be the leading coefficient of $N$ w.r.t $y_3$. Then
$A$ is either $-Z_2x_1+Z_1x_2$ or $-Z_2x_1$ or $Z_1x_2$. In any case, it is clear that
there exists an open subset of $U_0$, say $U_1$, such that for $\bar\tau\in U_1$,
$A(\bar\tau,y_1,y_2)$ does not vanish.
\item[(iii)]  Let $T(\bar x)=\operatorname{Dis}_{d}(g(\bar x,d))$. Note that $g$ is square-free
and primitive w.r.t $d$, and hence $T$ is not identically zero. Let $\Psi_2$ be the
curve defined by $T$ in $\mathbb C^2$ if $T$ is not constant and $\Psi_2=\emptyset$
otherwise. Then we consider the open subset $U_2=U_1\setminus\Psi_2$.

Now, let $\bar\tau\in U_2$. Then $g(\bar\tau,d)$ has exactly $\delta_d$ roots because of
$U_1$, being all different because of $U_2$. Proposition \ref{Prop03} implies that these
roots can be grouped in pairs, with elements in each pair differing only by
multiplication by $-1$. Let $\Theta(\bar\tau)=\{d_1,\ldots,d_{\mu}\}$ be a collection of
$\mu$ roots of $g(\bar\tau,d)$ where each $d_i$ is from one of these pairs. Also,
observe that because of $U_0$, $d_i\neq 0, \forall i=1\ldots,\mu$.

\item[(iv)] Now, let $\Upsilon$ be the set in Theorem \ref{Prop01}. Also consider the finite
(possibly empty) set $\tilde\Upsilon$ of values of $d$ such that for
$d_0\in\tilde\Upsilon$, $\mathcal O_{d_0}(\mathcal C)$ has a special component (see
section 5 in \cite{SS99}). Let
$\Psi_3=\cup_{d_0\in(\Upsilon\cup\tilde\Upsilon)}{\mathcal O}_{d_0}(\mathcal C)$, and
take $U_3=U_2\setminus\Psi_3$.

\item[(v)] Recall that ${\mathcal O}_{d}(\mathcal C)=
\left(\pi\left(\Omega(d)\right)\right)^*$. Let ${\mathcal M}={\mathcal O}_{d}(\mathcal
C)\setminus\pi\left(\Omega(d)\right)$. Let us see that, if $\mathcal M\neq\emptyset$,
then $\dim(\mathcal M)\leq 1$. For this purpose, let:
\[\Omega(d)=\Gamma_1\cup\cdots\cup\Gamma_s\]
where $\Gamma_i$ are the irreducible components of $\Omega(d)$. Let ${\mathcal
O}_{i}=(\pi(\Gamma_i))^*$. Then, since\vspace{2pt}

$
\mathcal M={\mathcal O}_{d}(\mathcal C)\setminus\pi(\Omega(d))=\\
=\left(\pi(\bigcup_{i=1}^s\Gamma_i)\right)^*\setminus\pi(\bigcup_{i=1}^s\Gamma_i)
=\bigcup_{i=1}^s\left(\pi(\Gamma_i)\right)^*\setminus\bigcup_{i=1}^s\pi(\Gamma_i)=\\
=\bigcup_{i=1}^s\left(\pi(\Gamma_i)^*\right)\setminus\bigcup_{i=1}^s\left(\pi(\Gamma_i)\right)
\subset\bigcup_{i=1}^s\left(\pi(\Gamma_i)^*\setminus\pi(\Gamma_i)\right) $\\

\noindent if $\dim(\mathcal M)>1$, there exists $i\in\{1,\ldots,s\}$ such that
$\dim\left(\pi(\Gamma_i)^*\setminus\pi(\Gamma_i)\right)>1$. Consider now the rational
map $\pi:\Gamma_i\mapsto\pi(\Gamma_i)^*$; note that both closed sets are irreducible. By
Theorem 7(ii) in \cite{Shafarevich77}, page 76. there exists a non-empty open subset $U$
of $\pi(\Gamma_i)^*$ such that the dimension of the fiber is invariant. Hence $\mathcal
M\subset\pi(\Gamma_i)^*\setminus U$, which is a contradiction, because
$\dim\left(\pi(\Gamma_i)^*\setminus U\right)\leq 1$. Now, we consider the projection
\[\pi_{\bar x}:\mathbb C^3\rightarrow\mathbb C^2; (x,d)\mapsto\bar x\]
Then $\Psi_4=(\pi_{\bar x}(\mathcal M))^*$ is either empty or $\dim(\Psi_4)\leq 1$. Let
us define $U_4=U_3\setminus\Psi_4$.

\item[(vi)]Consider the following resultants:
\[
R_i(\bar x)=\operatorname{Res}_{d}\left(g(\bar x,d),\dfrac{\partial g}{\partial
x_i}(\bar x,d)\right)
\]
for $i=1,2$. Note that $\dfrac{\partial g}{\partial x_i}$ cannot be identically zero,
because $\mathcal C$ is not a line. Also observe that $R_i$ cannot be identically zero,
since this would imply that $\dfrac{\partial g}{\partial x_i}(\bar x,d)$ and $g(\bar
x,d)$ have a common factor of positive degree in $d$. This factor cannot depend only on
$d$ because of the definition of the generic offset equation. Thus, this would imply
that for $d\not\in\Upsilon$ (the set in Theorem \ref{Prop01}), the offset has infinitely
many ramification points, and this is impossible since the offset cannot have multiple
components, and it cannot be a line because $\mathcal C$ is not a line. Let $\Phi_i$ be
the zero set of $R_i(\bar x)$ in $\mathbb C^2$. Take
$U_5=U_4\setminus(\Phi_1\cap\Phi_2)$. Now, if $\bar\tau\in U_5$, and $g(\bar\tau,d_0)=0$
with $d_0\not\in\Upsilon$, it follows that $\bar\tau$ is a regular point of ${\mathcal
O}_{d_0}(\mathcal C)$. Otherwise one has
\[g(\bar\tau,d_0)=\dfrac{\partial g}{\partial x_i}(\bar\tau,d_0)=0\]
for $i=1,2.$ This means that $R_i(\bar\tau)=0$ for $i=1,2$, contradicting the
construction of $U_5$.

\item[(vii)] Let $\{\tilde p_1,\ldots,\tilde p_r\}$ be the isotropic affine and regular points
of $\mathcal C$. This is a finite set because $\mathcal C$ is irreducible. For
$i=1,\ldots,r$, let $\gamma_i$ be the normal line to $\mathcal C$ at $\tilde p_i$. Let
$U_6=U_5\setminus\bigcup_{i=1}^r\gamma_i$.

\item[(viii)] If $(0:0:1)\in\mathcal C$ and it is regular, let $\Psi$ be the zero set in $\mathbb C^2$ of
\[N(0:0:1,\bar x)=-f_2(0,0)x_1+f_1(0,0)x_2\]
and define $U_7=U_6\setminus\Psi$; note that, since $(0:0:1)$ is regular in $\mathcal
C$.

\end{enumerate}
Let us see that $U=U_7$ satisfies the requirements. Let $\bar\tau\in U$, and let
$d_i\in\Theta(\bar\tau)$ (see the construction of $U_2$). Then $g(\bar\tau,\pm d_i)=0$.
Thus $(\bar\tau,\pm d_i)\in {\mathcal O}_{d}(\mathcal C)$ because of $U_3$. Moreover,
because of $U_4$, $\bar\tau\not\in\pi_{\bar x}(({\mathcal M}))^*$. Hence, $(\bar\tau,\pm
d_i)\in\pi(\Omega(d))$. Thus, there exist $\hat p_i\in\mathcal C$ and $u_0\in\mathbb C$
such that $(\hat p_i,\bar\tau,u_0)$ is a solution of $\mathfrak S_1(\pm d_i)$. In
particular, this implies that $\hat p_i$ is a solution of $\mathfrak S_3(\bar\tau)$, and
that $\hat p_i$ generates $\bar\tau$ in ${\mathcal O}_{\pm d_i}(\mathcal C)$. Let
$\hat\Gamma=\{\hat p_1,\ldots,\hat p_{\mu}\}$. Observe that $\hat p_i\in\mathcal C$ and
it is affine. Moreover, since $(\hat p_i,\bar\tau,u_0)$ is a solution of $\mathfrak
S_1(\pm d_i)$, then $\hat p_i$ is non-isotropic on $\mathcal C$.  Now, since $d_i\neq
d_j$ for $i\neq j$ (see the construction of $U_2$), and since $\hat p_i$ belongs to a
circle of radius $d_i$ and centered at $\bar\tau$, one concludes that $\hat p_i\neq\hat
p_j$. So, statement (1) and (4) hold. Statement (3) follows from the construction of
$U_2$.

The existence part of Statement (2) follows from the construction of $U_6$ and $U_7$. It
remains only to prove that, for $\bar\tau\in U$, $\hat\Gamma(\bar\tau)$ contains {\em
all} the affine and non-isotropic solutions of $\mathfrak S_3(\bar\tau)$. Suppose that
$\tilde p$ is an affine non-isotropic point of $\mathcal C$ such that $N(\tilde
p,\bar\tau)=0$ and $\tilde p\not\in\hat\Gamma(\bar\tau)$. Because of $U_2$, it follows
that $\tilde p$ generates $\bar\tau\in{\mathcal O}_{\pm d_i}(\mathcal C)$ for some
$d_i\in\Theta(\bar\tau)$. Then, we could take places of $\mathcal C$ at both $\tilde p$
and $\hat p_i$ and lift them to places of the offset at $\bar\tau$. Since ${\mathcal
O}_{\pm d_i}(\mathcal C)$ has no special component, these two places cannot lift to the
same place of the offset. But if they lift to different places, it follows that
$\bar\tau$ is not regular in ${\mathcal O}_{\pm d_i}(\mathcal C)$, and this contradicts
the construction with $U_5$. \hfill\qed}

In the next definition we extend the terminology of fake and non-fake points to this
degree problem.
\begin{defn}
Let $U$ be as in Theorem \ref{Theorem06}. We denote:
\[d{\mathcal F}=\bigcap_{\bar \tau\in U}\left[\overline{{\mathcal N}(\bar
\tau)} \cap \overline{{\mathcal C}}\right]
\]
The points of the set $d{\mathcal F}$ are called {\sf $d$-fake points}. For $\bar
\tau\in U$, the points in $\left(\overline{{\mathcal N}(\bar \tau)} \cap
\overline{{\mathcal C}}\right)\setminus d{\mathcal F}$ are called {\sf non $d$-fake
points}.
\end{defn}
The next step in the strategy consists in showing the invariance of the set of $d$-fake
points. This is established in the next proposition (compare to Proposition
\ref{CaracterizacionInvarianteFakes}).
\begin{prop}[\sf Invariance of the $d$-fake points]\label{Lemma02a}
Let $U$ be as in Theorem \ref{Theorem06}. The set $d\mathcal F$ is finite. Moreover,
\[d\mathcal F=\operatorname{Sing}_a(\overline{\mathcal C})
\cup\operatorname{Iso}_{\infty}(\overline{\mathcal C})
\]
where $\operatorname{Sing}_a(\overline{\mathcal C})$ is the affine singular locus of
$\overline{\mathcal C}$ and $\operatorname{Iso}_{\infty}(\overline{\mathcal C})$ is the
set of isotropic points at infinity of $\overline{\mathcal C}$; that is, the set of
points of $\mathcal C$ that satisfy $y_3=0$ and $F_1^2+F_2^2=0$.
\end{prop}
{\bf Proof. }{\small

Let $p=(a:b:c)\in d\mathcal F$. Then $p\in\overline{\mathcal C}$ and $N(p,\bar\tau)=0$
for every $\bar\tau\in U$. Thus, considering $N(\bar y_H,\bar x)\in\mathbb C[\bar
y_H][\bar x]$, one has that:
\[-F_2(p)c=0,\quad F_1(p)c=0,\quad F_2(p)a-F_1(p)b=0.\]
If $c\neq 0$, then $p$ is affine and $F_1(p)=F_2(p)=0$. Thus
$p\in\operatorname{Sing}_a(\overline{\mathcal C})$. If $c=0$, then using Euler's
identity
\[F_1(p)a+F_2(p)b=\deg(F)F(p)=0\]
From this relation and $F_2(p)a-F_1(p)b=0$ one has that $F_1^2(p)+F_2^2(p)=0$. Thus
$p\in\operatorname{Iso}_{\infty}(\overline{\mathcal C})$. Therefore $d\mathcal F\subset
\operatorname{Sing}_a(\overline{\mathcal
C})\cup\operatorname{Iso}_{\infty}(\overline{\mathcal C})$.

Conversely, let $p\in\operatorname{Sing}_a(\overline{\mathcal
C})\cup\operatorname{Iso}_{\infty}(\overline{\mathcal C})$. If
$p=(a:b:c)\in\operatorname{Sing}_a(\overline{\mathcal C})$, then $p\in\overline{\mathcal
C}$ and for every $\bar\tau\in U$ one has $N(p,\bar\tau)=0$. Thus, $p\in d\mathcal F$.
If $p\in\operatorname{Iso}_{\infty}(\overline{\mathcal C})$, then $p\in\mathcal C$,
$c=0$, and $F_1^2(p)+F_2^2(p)=0$. Using Euler's identity as before one has
$F_1(p)a+F_2(p)b=0$. From these relations one gets $N(p,\bar\tau)=F_2(p)a-F_1(p)b=0$ for
all $\bar\tau\in U$. Thus, $p\in d\mathcal F$.

The finiteness of $d\mathcal F$ follows from the equality $d\mathcal F=
\operatorname{Sing}_a(\overline{\mathcal
C})\cup\operatorname{Iso}_{\infty}(\overline{\mathcal C})$.}

\begin{rem}\label{remIsotropiadFakesInfinito}
\begin{enumerate}
\item[]
\item The proof of Proposition \ref{Lemma02a} shows that if $p$ is a point at infinity of
$\mathcal C$, and for some $\bar\tau\in U$, $p\in\overline{{\mathcal N}(\bar
\tau)}\cap\overline{{\mathcal C}}$, then
$p\in\operatorname{Iso}_{\infty}(\overline{\mathcal C})$
\item From the definition of $d\mathcal F$ it follows that for any non empty open subset
$\tilde U\subset U$, one has
\[d{\mathcal F}=
 \bigcap_{\bar \tau\in U}\left[\overline{{\mathcal N}(\bar \tau)} \cap
\overline{{\mathcal C}}\right]=
 \bigcap_{\bar \tau\in\tilde U}\left[\overline{{\mathcal
N}(\bar \tau)} \cap \overline{{\mathcal C}}\right]
\]
\end{enumerate}
\end{rem}

\begin{prop}[\sf Characterization of the $d$-fake points]\label{Caracterizacion_d_Fakes}
Let $U$ be as in Theorem \ref{Theorem06}. With the notation of Theorem \ref{Theorem06},
for each $\bar\tau\in U$, it holds that
\begin{enumerate}
\item $\overline{{\mathcal N}(\bar \tau)}\cap\overline{{\mathcal C}}=\hat\Gamma(\bar\tau)\cup
 d{\mathcal F} $
\item $\hat\Gamma(\bar\tau)\cap  d{\mathcal F}=\emptyset$
\end{enumerate}
\end{prop}
{\bf Proof. }{\small

Let $\bar\tau\in U$.

(1) Let $p=(a:b:c)\in\overline{{\mathcal N}(\bar \tau)}\cap\overline{{\mathcal C}}$. If
$c=0$, then by Remark \ref{remIsotropiadFakesInfinito} (1), one has
$p\in\operatorname{Iso}_{\infty}(\overline{\mathcal C})$, and by Proposition
\ref{Lemma02a}, $p\in d\mathcal F$. If $c\neq 0$ and $p\in
\operatorname{Sing}_a(\overline{\mathcal C})$, then again by Proposition \ref{Lemma02a},
$p\in d\mathcal F$. If $c\neq 0$ and $p\not\in \operatorname{Sing}_a(\overline{\mathcal
C})$, then $p$ is an affine regular point of $\mathcal C$. By Theorem \ref{Theorem06},
then $p\in\hat\Gamma(\bar\tau)$. Thus, in any case, $p\in d\mathcal
F\cap\hat\Gamma(\bar\tau)$. The reverse inclusion is trivial.

(2) This follows from Proposition \ref{Lemma02a}. \hfill\qed}

\begin{rem}\label{remLosNodFakesSonGamma}
Proposition \ref{Caracterizacion_d_Fakes} shows that if $\bar\tau\in U$, then the set of
non $d$-fake points is precisely $\hat\Gamma(\bar\tau)$. In particular,
\[\operatorname{Card}([\overline{{\mathcal N}(\bar\tau)}\cap \overline{{\mathcal
C}}]\setminus d{\mathcal
F})=\operatorname{Card}(\hat\Gamma(\bar\tau))=\mu=\dfrac{\delta_d}{2}\]
\end{rem}

The next proposition gathers the information we need when applying Bezout's Theorem to
the curves $\overline{\mathcal C}$ and $\overline{\mathcal N}(\bar\tau)$ (compare to
Proposition \ref{Prop02}).

\begin{prop}\label{Prop04}
There exists a non-empty open subset $\tilde U\subset U$, where $U$ is as in Theorem
\ref{Theorem06}, such that for every $\bar\tau\in\tilde U$ the following hold:
\begin{enumerate}
\item  $\deg({\mathcal N}(\bar y_H,\bar x))$ does not depend on $\bar x$,
\item ${\mathcal C}$ and ${\mathcal N}(\bar\tau)$ have no common component,
\item  if $\hat p$ is a non $d$-fake point,
then $\operatorname{mult}_{\hat p}(\overline{\mathcal C},\overline{\mathcal
N(\bar\tau)})=1$.
\item Let $N(\bar y_H,\bar x)$ be considered as an element of $(\mathbb
C[\bar y_H])[\bar x]$:
\[N(\bar y_H,\bar x)=Z_{1,0}(\bar y_H)x_1+Z_{0,1}(\bar y_H)x_2+Z_{0,0}(\bar y_H)\]
where:
\[
\begin{cases}
Z_{1,0}(\bar y_H)=-F_{2}y_3\\
Z_{0,1}(\bar y_H)=F_1y_3,\\
Z_{0,0}(\bar y_H)=F_2y_1-F_1y_2,
\end{cases}
\]
and let $\mathcal J_{\alpha}$ be the zero in $\mathbb C^2$ set of $Z_{\alpha}(\bar
y_H)$. Then it holds that:
\[\bigcap_{\alpha}(\overline{{\mathcal C}} \cap \overline{{\mathcal J_{\alpha}}})\subset
d{\mathcal F}.\]
\item $(0:0:1)\not\in\left(\overline{{\mathcal
N(\bar\tau)}} \cap \overline{{\mathcal C}}\right)\setminus{d\mathcal F}$
\end{enumerate}
\end{prop}

{\bf Proof. }{\small
\begin{enumerate}
\item See step (ii) in the proof of Theorem \ref{Theorem06}.
\item Let us consider $n$ as a polynomial in ${\mathbb C}[y_1,y_2][x_1,x_2]$. If $n$ and
$f$ have a common factor, one has that $f_1=f_2=0$ for every point of $\mathcal C$,
which is a contradiction since $\mathcal C$ is irreducible.
\item Let
\[P(t)=(y_1(t),y_2(t))\]
with
$$\begin{cases}
y_1=a_0+a_1t+\cdots\\
y_2=b_0+b_1t+\cdots
\end{cases}$$
be a place of ${\mathcal C}$ centered at $\hat p$. Then the multiplicity of intersection
$\operatorname{mult}_{\hat p}(\overline{\mathcal C},\overline{\mathcal N(\bar\tau)})$ is
equal to the order of $n(P(t),\bar\tau)$. Let now
\[
\begin{cases}
f_1(P(t))=\alpha_0+\alpha_1 t+\cdots\\
f_2(P(t))=\beta_0+\beta_1 t+\cdots
\end{cases}
\]
Note that $\alpha _{0}^{2}+\beta _{0}^{2}\neq 0$ because $\hat p$ is non $d$-fake.
Besides, since the point $\bar\tau=(\tau_1,\tau_2)$ is generated by $\hat p$ in
${\mathcal O}_{d_i}(\mathcal C)$, one has:
\[\tau _{1}=a_{0}+d_i\frac{\alpha _{0}}{\sqrt{\left( \alpha _{0}^{2}+\beta _{0}^{2}\right) }}\]
\[\tau _{2}=b_{0}+d_i\frac{\beta _{0}}{\sqrt{\left( \alpha _{0}^{2}+\beta _{0}^{2}\right) }}\]
Substituting the above expressions in $n$ we arrive at:
$$
n(P(t),\bar\tau)=\left( -\alpha _{0}b_{1}+\alpha _{1}d_i\frac{\beta _{0}}{\sqrt{\left(
\alpha _{0}^{2}+\beta _{0}^{2}\right) }}+\beta _{0}a_{1}-\beta _{1}d_i\frac{\alpha
_{0}}{\sqrt{\left( \alpha _{0}^{2}+\beta _{0}^{2}\right) }}\right) t+\cdots
$$
(the order zero term vanishes identically.) Now, we will suppose that we have
$\operatorname{mult}_{\hat p}({\mathcal C},{\mathcal N})>1$ and we will arrive at a
contradiction.  This would imply that
$$ -\alpha _{0}b_{1}+\alpha _{1}d_i\frac{\beta _{0}}{\sqrt{\left( \alpha
_{0}^{2}+\beta _{0}^{2}\right) }}+\beta _{0}a_{1}-\beta _{1}d_i\frac{\alpha
_{0}}{\sqrt{\left( \alpha _{0}^{2}+\beta _{0}^{2}\right) }}=0
$$
From this one gets:
\[
(-\alpha _{1}\beta _{0}+\beta _{1}\alpha _{0})d_i=-\sqrt{\left( \alpha _{0}^{2}+\beta
_{0}^{2}\right) }(\alpha _{0}b_{1}-\beta _{0}a_{1})
\]
Now observe that $\left( a_{1},b_{1}\right)$ is a tangent vector to $\mathcal C$ at
$\hat p$, and $\left( \alpha _{0},\beta _{0}\right) $ is a normal at the same point.
Thus $a_{1}\alpha _{0}+b_{1}\beta _{0}=0$. Thus, if $-\alpha _{1}\beta _{0}+\beta
_{1}\alpha _{0}=0$, since $\alpha _{0}^{2}+\beta _{0}^{2}\neq 0$, one obtains:
$$\begin{cases}
\alpha _{0}b_{1}-\beta _{0} a_{1}=0\\
\beta _{0}b_{1}+\alpha _{0} a_{1}=0
\end{cases}$$
It follows that $a_1=b_1=0$, which is a contradiction, since $\hat p$ is regular in
$\mathcal C$. Thus, we have shown that $-\alpha _{1}\beta _{0}+\beta _{1}\alpha _{0}\neq
0$.  And therefore
\[
d_i=\sqrt{\left( \alpha _{0}^{2}+\beta _{0}^{2}\right) }\dfrac{\alpha _{0}b_{1}-\beta
_{0}a_{1}}{\alpha _{1}\beta _{0}-\beta _{1}\alpha _{0}}
\]
Now, as in the proof of Proposition \ref{Prop03}, we can offset the place $P(t)$ to get
a place of  ${\mathcal O}_{d_i}(\mathcal C)$ centered at $\bar\tau$
$$O(t)=(O_1(t),O_2(t))=
\left(
 y_{1}(t)+d_i\frac{f_{1}(t)}{\sqrt{f_{1}^{2}(t)+f_{2}^{2}(t)}},
 y_{2}(t)+d_i\frac{f_{2}(t)}{\sqrt{f_{1}^{2}(t)+f_{2}^{2}(t)}}
\right)
$$
Substituting the above expressions by $y_1(t), y_2(t), f_1(t), f_2(t)$ and $d_i$ one
has, after simplifying the expression:
\[
O_1(t)=\tau_1+ \left( a_{1}\alpha _{0}+b_{1}\beta _{0}\right) \frac{\alpha _{0}}{\alpha
_{0}^{2}+\beta _{0}^{2}}  t+\cdots
\]
Similarly
\[
O_2(t)=\tau_2+\allowbreak \left( a_{1}\alpha _{0}+b_{1}\beta _{0}\right) \frac{\beta
_{0}}{\alpha _{0}^{2}+\beta _{0}^{2}}t+\cdots
\]
Since $a_{1}\alpha _{0}+b_{1}\beta _{0}=0$, this would imply that $\bar\tau$ is not
regular in ${\mathcal O}_{d_i}(\mathcal C)$, contradicting the construction of the open
set $U$.
\item If $p=(a:b:c)\in\bigcap_{\alpha}(\overline{{\mathcal C}} \cap \overline{{\mathcal
J_{\alpha}}})$, then $cF_1(p)=0$ and $cF_2(p)=0$. If $c\neq 0$, it follows that
$F_1(p)=F_2(p)=0$. If $c=0$, $F_1^2(p)+F_2^2(p)=0$ follows by Remark
\ref{remIsotropiadFakesInfinito} (1). In either case, by Proposition \ref{Lemma02a},
$p\in d\mathcal F$.

\item This follows from statement (1) in Theorem \ref{Theorem06}. \hfill\qed
\end{enumerate}}

\section{Degree formulae for the distance}\label{secDistanceFormula}

As a consequence of the results in the previous section, we derive the following formula
for computing $\delta_d$.
\begin{thm}[\sf Degree formula for the distance]\label{Theorem05}
\[\delta_d=\deg_{d}({\mathcal O}_d({\mathcal C}))=2\,
\deg_{\{y_1,y_2\}}\left(\operatorname{PP}_{\{x_1,x_2\}}\left(\operatorname{Res}_{y_3}(F(\bar
y_H),N(\bar y_H,\bar x))\right)\right)\] We recall that $F$ is the homogeneous implicit
equation of the curve, and $N$ is the polynomial introduced after Remark \ref{RemarkMu}.
\end{thm}
{\bf Proof of Theorem.}{\small

In order to prove the theorem, we apply Theorem \ref{Theorem07}. Let $\mathcal
D=\mathcal C$, $Z(\bar y_H,\bar u)=N(\bar y_H,\bar x)$, where $\bar x=(x_1,x_2)$, and
$\Xi=\tilde U$, where $\tilde U$ is as in Proposition \ref{Prop04}. We check that all
the hypothesis are satisfied:
\begin{itemize}
\item $\mathcal C$ is irreducible and it is not a line by assumption.
\item $N$ can be written as
\[N=(-F_2x_1+F_1x_2)y_3+(y_1F_2-y_2F_1)\]
Thus, since $F_1$ and $F_2$ are not identically zero, $S$ depends on $y_3$.
\item (1) and (2) in Theorem \ref{Theorem07} follow from (1) and (2)  in
Proposition \ref{Prop04}.
\item The equality $d{\mathcal F}=\bigcap_{\bar x\in\tilde U} \left[\overline{{\mathcal N}(\bar x)} \cap \overline{{\mathcal
C}}\right]$ follows from  Remark \ref{remIsotropiadFakesInfinito}(2).
\item In this situation, hypothesis (3), (4) and (5) in Theorem \ref{Theorem07} follows from
Proposition \ref{Prop04} (3), (4) and (5).
\end{itemize}
Then, Theorem \ref{Theorem07} implies that there exists a non-empty open
$U^*\subset\tilde U$ such that for $\bar\tau\in U^*$
\[\operatorname{Card}([\overline{{\mathcal N}(\bar\tau)}\cap \overline{{\mathcal C}}]\setminus {d\mathcal F})=
\deg_{\{y_1,y_2\}}\left(\operatorname{PP}_{\bar x}\left(\operatorname{Res}_{y_3}(F(\bar
y_H),N(\bar y_H,\bar x))\right)\right)
\]
Now the theorem follows from Remark \ref{remLosNodFakesSonGamma}. \hfill\qed}

\appendix{\bf Appendix: Table of Offset degrees}\label{Apenddix}

In the following table we list, for some curves, the total degree $\delta$ w.r.t
$\{x_1,x_2\}$ of the generic offset equation $g(x_1,x_2,d)$, its partial degrees
$\delta_1$ and $\delta_2$ w.r.t $x_1$ and $x_2$, respectively, and the degree $\delta_d$
w.r.t $d$.

\begin{tabular}{|c||c|c|c|c|c|}
\hline
Curve $\mathcal C$&Equation $f(y_1,y_2)=0$&$\delta$&$\delta_1$&$\delta_2$&$\delta_d$\\
\hline\hline
 Circle&$y_1^2+y_2^2-r^2=0$     &   4   &   4    &   4    &    4   \\
\hline
 Parabola&$y_2+a+by_1+cy_1^2=0$     &   6   &   6    &   4    &    6   \\
\hline
 Ellipse&$y_1^2/a^2+y_2^2/b^2-1=0$     &   8   &   8    &   8    &    8   \\
\hline
 Hyperbola&$y_1^2/a^2-y_2^2/b^2-1=0$     &   8   &   8    &   8    &    8   \\
\hline
 Hyperbola&$y_1y_2-1=0$     &   8   &   6    &   6    &    8   \\
\hline
 Cubic Cusp&$y_1^3-y_2^2=0$     &   8   &   8    &   6    &    8   \\
\hline
 Folium&$y_1^3+y_2^3-3y_1y_2=0$&14&14&14&14\\
\hline
 Conchoid&$(y_1-1)(y_1^2+y-2^2)+y_1^2=0$&8&8&6&8\\
\hline
 A cubic&$y_1^3+y_2^3-y_1y_1-1=0$&18&18&18&18\\
\hline
 Epitrochoid&$y_2^4+2y_1^2y_2^2-34y_2^2+y_1^4-34y_1^2+96y_1-63=0$&10&10&10&8\\
\hline
 Cardioid&$(y_1^2+4y_2+y_2^2)^2-16y_1^2-16y_2^2=0$&8&8&8&6\\
\hline
 Rose (three petals)&$(y_1^2+y_2^2)^2+y_1(3y_2^2-y_1^2)=0$&14&14&12&12\\
\hline
 Ramphoid Cusp&$y_1^4+y_1^2y_2^2-2y_1^2y_2-y_1y_2^2+y_2^2=0$&14&14&10&14\\
\hline
 Lemniscate&$(y_1^2+y_2^2)^2-2(y_1^2-y_2^2)=0$&12&12&12&12\\
\hline
 Scarabeus&$(y_1^2+y_2^2)(y_1^2+y_2^2+y_1)^2-(y_1^2-y_2^2)^2=0$&18&18&18&14\\
\hline
\end{tabular}

\end{document}